\documentclass[12pt]{article}

\usepackage{natbib}
\usepackage{verbatim,color,amssymb}
\usepackage{enumerate}
\usepackage{amsmath}					
\usepackage{multirow}
\usepackage{setspace}

\usepackage{graphicx}
\usepackage{color}

\newtheorem{Thm}{\underline{\bf Theorem}}

\newtheorem{Cond}{\underline{\bf Conditions}}
\newtheorem{Proof}{Proof}

\newtheorem{Lem}{\underline{\bf Lemma}}

\newcommand{\MyProof}{\noindent\textbf{Proof. }}

\def\rR{\mathbb{R}}

\def\C{{\cal C}}

\def\F{{\cal F}}
\def\G{{\cal G}}
\def\K{{\cal K}}
\def\calP{{\cal P}}
\def\S{{\cal S}}

\def\calP{{\cal P}}

\def\wh{\widetilde}

\def\wh{\widehat}

\def\th{^{th}}
\def\Pr{\hbox{Pr}}

\def\Beta{\hbox{Beta}}
\def\Cauchy{\hbox{Cauchy}}

\def\Exp{\hbox{Exp}}

\def\IG{\hbox{Inv-Ga}}

\def\Normal{\hbox{Normal}}

\def\Wbll{\hbox{Weibull}}

\def\ANNALS{{\it Annals of Statistics}}

\def\BIOK{{\it Biometrika}}

\def\CANADAJS{{\it Canadian Journal of Statistics}}

\def\EJS{{\it Electronic Journal of Statistics}}

\def\JASA{{\it Journal of the American Statistical Association}}

\def\P_25_ICML{{\it Proceedings of the 25th international conference on Machine learning}}

\def\STATL{{\it Statistics and Probability Letters}}

\def\SSNC{{\it Statistica Sinica}}

\def\refhg{\hangindent=20pt\hangafter=1}
\def\refmark{\par\vskip 2mm\noindent\refhg}

\def\var{\hbox{var}}

\def\refhg{\hangindent=20pt\hangafter=1}
\def\refmark{\par\vskip 2mm\noindent\refhg}

\def\bse{\begin{eqnarray*}}
\def\ese{\end{eqnarray*}}
\def\be{\begin{eqnarray}}
\def\ee{\end{eqnarray}}
\def\bq{\begin{equation}}
\def\eq{\end{equation}}

\def\wh{\widehat}
\def\wt{\widetilde}

\def\b1e{{\mathbf e}}

\def\bW{{\mathbf W}}

\def\bX{{\mathbf X}}

\newcommand{\abs}[1]{\left\vert#1\right\vert}

\newcommand{\norm}[1]{\left\Vert#1\right\Vert}

\DeclareMathSymbol{\precapprox}{\mathrel}{AMSb}{119}

\def\authorfootnote#1{{\let\thefootnote\relax\footnotetext{#1}}}

\allowdisplaybreaks

\begin{document}
\thispagestyle{empty}
\baselineskip=25pt

\begin{center}
{\Large{\bf Adaptive Posterior Convergence Rates in Bayesian Density Deconvolution with Supersmooth Errors}}
\end{center}

\baselineskip=13pt

\vskip 2mm
\begin{center}
Abhra Sarkar\\
Department of Statistics, Texas A\&M University, 3143 TAMU, College Station, TX 77843-3143 USA\\
abhra@stat.tamu.edu\\
\hskip 5mm \\
Debdeep Pati\\
Department of Statistics, Florida State University, Tallahassee, FL\\ 32306-4330 USA\\
debdeep@stat.fsu.edu\\
\hskip 5mm \\
Bani K. Mallick\\
Department of Statistics, Texas A\&M University, 3143 TAMU, College Station, TX 77843-3143 USA\\
bmallick@stat.tamu.edu\\
\hskip 5mm \\
Raymond J. Carroll\\
Department of Statistics, Texas A\&M University, 3143 TAMU, College Station, TX 77843-3143 USA\\
carroll@stat.tamu.edu\\
\end{center}

\vskip 2cm
\begin{center}
{\Large{\bf Abstract}}
\end{center}
Bayesian density deconvolution using nonparametric prior distributions is a useful  alternative to the frequentist kernel based deconvolution estimators due to its potentially wide range of applicability, straightforward uncertainty quantification and generalizability to more sophisticated models. This article is the first substantive effort to theoretically quantify the behavior of the posterior in this recent line of research.  In particular, assuming a known supersmooth error density, a Dirichlet process mixture of Normals on the true density leads to a posterior convergence rate same as the minimax rate $(\log n)^{-\eta/\beta}$ adaptively over the smoothness $\eta$ of an appropriate H\"{o}lder space of densities, where $\beta$ is the degree of smoothness of the error distribution.  Our main contribution is achieving adaptive minimax rates with respect to the $L_p$ norm for $2 \leq p \leq  \infty$ under mild regularity conditions on the true density. En route, we develop tight concentration bounds for a class of kernel based deconvolution estimators which might be of independent interest.  
\baselineskip=12pt

\baselineskip=12pt
\par\vfill\noindent
\underline{\bf Some Key Words}: Density deconvolution, Dirichlet process mixture models, Measurement errors, Posterior convergence rates.
\par\medskip\noindent
\underline{\bf Short Title}: Posterior Convergence Rates in Density Deconvolution

\clearpage\pagebreak\newpage
\pagenumbering{arabic}
\newlength{\gnat}
\setlength{\gnat}{16pt}
\baselineskip=\gnat

\section{Introduction}
A density deconvolution problem is a specialized density estimation ($f_{X}$ of a random variable $X$) when precise observations on $X$ are not available, 
but observations on $W$, a contaminated proxy for $X$, contaminated with additive measurement error $U$, are available.
The data generating model is thus given by
\be  \label{eq: deconv}
W_{i} =  X_{i} + U_{i}, \quad i=1, \ldots, n.
\ee
Assuming $X$ and $U$ to be independent, in terms of densities, the observations are generated from the convolution $f_{W}(w) = (f_{X}\star f_{U}) (w) = \int f_{X}(w-u)f_{U}(u)du$, where $f_{U}$ and $f_{W}$ denote the densities of $U$ and $W$, respectively. 
In this article we assume that error density $f_{U}$ is completely known.

To solve the deconvolution problem in a Bayesian nonparametric framework, 
a prior distribution, denoted here by the generic notation $\Pi$, is assigned to the unknown density of interest $f_{X}$. 
Given a random sample $\bW_{1:n} = \{W_{1},\dots,W_{n}\}$  from $f_{W}$, 
Bayesian inference is then based on the posterior distribution is given by 
\be \label{eq: posterior}
\Pi(B\mid \bW_{1:n}) = \frac{\int_{B}\prod_{i=1}^{n}\int f_{X}(W_{i}-U_{i})f_{U}(U_{i}) dU_{i} d\Pi(f_{X})}{\int \prod_{i=1}^{n}\int f_{X}(W_{i}-U_{i})f_{U}(U_{i})dU_{i} d\Pi(f_{X})}.
\ee
In this article we wish to study 
the consistency properties of the posterior as $n\to \infty$ 
under the frequentist assumption of an underlying true density $f_{0X}$ for $X$.

Frequentist deconvolution estimators, in particular deconvoluting kernel type estimators (DKE, see Section \ref{sec: DKE}), have been extensively researched in the literature.
Optimal point wise convergence rates have been studied by Carroll and Hall (1988), Stefanski and Carroll (1990) and Fan (1991b), among others.
Global convergence rates were studied in Fan (1988) for weighted $L_p$ norm.
The results of Fan (1988, 1991a, 1991b) show that the convergence rates depend on the smoothness of the error distribution: the smoother the error distribution, the more difficult it is to recover the density of interest.
In particular, when the true density belongs to a normed H\"{o}lder class of smoothness $\eta$ and 
the error distribution is supersmooth (the characteristic function has exponential decay) with smoothness $\beta$,
the minimax optimal rate of convergence is $(\log n)^{-\eta/\beta}$ and is attained by the DKE.

For density estimation problems, where in contrast there is no measurement error and accurate measurements on $X$ are available,
Bayesian nonparametric techniques including Dirichlet process mixture models (DPMM) (Ferguson, 1973; Lo, 1984; Escobar and West, 1995), 
where the unknown density is modeled as a mixture of normals with a Dirichlet process prior on the mixing distribution, 
have been hugely successful.
Flexibility and richness aside, 
the immense popularity of these methods can be attributed largely to the development of 
sophisticated computational machinery that has made implementation of these techniques routine in various applied problems.
To establish further credibility of such methods, frequentist consistency properties have also been given substantial attention in the literature 
and results of the type
\be\label{eq: rate general}
E_{f_{0X}} \Pi_{n}\{d(f_{0X},f_{X})>\xi_{n} \mid \bX_{1:n}\} \to 0
\ee
have been established, where $\bX_{1:n}$ denotes a set of precise measurements on $X$, $\xi_{n} \to 0$ and $d$ denotes a distance metric.
Such posterior convergence results imply the frequentist convergence rate $\xi_{n}$ for the associated Bayesian procedure.  
Posterior consistency and optimal rates of posterior convergence in the Hellinger metric have been studied by Ghosal, Ghosh and Ramamoorthi (1999), Ghosal, Ghosh and van der Vaart (2000), Ghosal and van der Vaart (2007a, 2007b), Shen and Wasserman (2001), among others.
More recently, Gin\'{e} and Nickl (2011) generalized the results to $L_{p}$ norms for $1\leq p \leq \infty$.

Bayesian nonparametric density estimation approaches, such as the DPMM,  can be readily adapted to the problem of density deconvolution, with practically no additional computational effort.
For recent methodological contributions in this direction see Sarkar, et al. (2013).
However, since in a deconvolution context the density of interest is different from the data generating density, theoretical investigation of consistency properties of the posterior is 
substantially different. In particular, it is not immediately clear whether the same formulation of the DPMM as in Ghosal and van der Vaart (2007) can lead to adaptive minimax optimal rates even in the case of density deconvolution.  

In this article, we show that when measurement errors are supersmooth, 
under some mild conditions on the true density,  the posterior obtained from a suitably chosen DPMM on $f_X$ converges to the truth at the minimax rate $(\log n)^{-\eta/\beta}$.  One of our main contributions is to formulate  the convergence of  (\ref{eq: posterior}) in the $L_{p}$ norm for $2 \leq p \leq \infty$. The set of sufficient conditions are milder compared to the case of usual density estimation in Ghosal and van der Vaart (2007) in that we only require polynomially decaying tails of the true density of the $X$. Moreover we achieve adaptivity to all smoothness levels of the true density of the $X$ based on realistic prior assumptions.  To the best of our knowledge,  achieving adaptive minimax rates with respect to $L_{\infty}$ norm is an open problem even in density estimation.

Since density deconvolution can be viewed as an inverse problem, 
our work is related to the recent works of Knapik, et al. (2011) and Ray (2013).  
While the work of Knapik, et al. (2011) is restricted to conjugate priors, Ray (2013) considers only periodic function deconvolution using wavelets.  
Although, we follow the general recipe given in Theorem 3.1 of Ray (2013) as the sufficient conditions for posterior convergence in an inverse problem,  
substantial technical hurdles remain.  
One of the main ingredients of Theorem 3.1 of Ray (2013) is to exploit the concentration properties of frequentist estimators to 
construct  test functions with type-I and type-II error bounds of the type $\exp(-Cn\epsilon_{n}^{2})$ for the testing problem
\be \label{eq: testing problem}
H_{0}: f_{X}=f_{0X} ~~~vs~~~ H_{A}: f_{X}\in \{f: d(f,f_{0X})>\xi_{n}\}.
\ee  
Ray (2013) used concentration properties of thresholded wavelet based estimators based on standard results on concentration of Gaussian priors. 
On the contrary, analogous results for DKE estimators suited to density deconvolution problems are lacking. 
One of our key technical contributions is to develop sharp concentration inequalities of the DKE to construct  tests for  (\ref{eq: testing problem}).

The article is organized as follows. 
Section \ref{sec: main thms} gives the main results.
A set of sufficient conditions, used to prove the main results, are provided in Section \ref{sec: suff thm}.
Section \ref{sec: aux results} details some auxiliary results used in the construction of the test function and in the proof of the main results.
Section \ref{sec: pf main thm} gives a proof of the main theorems combining the auxiliary results of Section \ref{sec: aux results}.
The optimal rate of convergence of deconvolution estimators for supersmooth errors is extremely slow.
In Section \ref{sec: acc rates} we discuss how small $\sigma$ should be for deconvolution with supersmooth errors to be practically feasible 
and for the deconvolution estimator to converge as fast as an ordinary density estimator.
Section \ref{sec: discussion} concludes the article with a discussion. 
The proof of the theorem stating sufficient conditions is provided in the Appendix.

\section{Notations}

Let $C(\mathbb{R})$ denote the space of all real valued bounded continuous functions on $\rR$.
For $2 \leq p \leq \infty$ and $\eta > 0$, let $C_{p}^{\eta}(\mathbb{R})$ denote the normed H\"{o}lder space comprising functions $f \in C(\mathbb{R})$ that have finite derivatives $f^{(k)}$ up to order $k\leq \lfloor \eta \rfloor$, 
with $\norm{f^{(\lfloor \eta \rfloor)}(x+y)-f^{( \lfloor \eta \rfloor)}(x)}_{p} \precsim  \abs{y}^{(\eta - \lfloor \eta \rfloor)}$ for any $x, y \in \rR$.
Denote the Fourier transform of a function $f$ by $\phi_{f}(t) =\int e^{i t x} f(x) dx$ and the inverse Fourier transform using the convention $f(x) = (2\pi)^{-1} \int e^{-i t x} \phi(t) dt$. 
A density $f$ is said to be supersmooth of order $\beta>0$ if, for some constants $\beta_{0}$ and $d_{0}, d_{1}>0$,
\bse \label{eq: supersmooth defn}
d_0 \abs{t}^{\beta_0} \exp (- \abs{t}^{\beta} / \varrho) \leq \abs{\phi_{f}(t)} \leq
d_1 \abs{t}^{\beta_1} \exp (- \abs{t}^{\beta} / \varrho), \, \text{as} \, t \to \infty. 
\ese
We assume $U = \sigma \wt{U}$, where $\sigma$ is a scale parameter parametrizing $f_{U}$ and $\wt{U}$ denotes a scaled version of the measurement errors.
For notational convenience $f_{U}$ will henceforth be denoted by $\psi_{\sigma}$, 
and the corresponding characteristic function (CF) will be denoted by $\phi_{\sigma}$.  
When $\sigma = 1$, $\psi_{\sigma}$ and $\phi_{\sigma}$ will be denoted simply by $\psi$ and $\phi$, respectively.  
Throughout we assume that $\psi_{\sigma}$ is supersmooth of order $\beta$. 
The CFs of $f_{0X}$ and $f_{0W}$ will be denoted simply by $\phi_{0X}$ and $\phi_{0W}$, respectively.
Throughout the article $\Normal(0,\sigma^{2})$ and $\Cauchy(0,\sigma^{2})$ denote a Normal and a Cauchy density, respectively, with median $0$ and scale $\sigma$. 
For $a>0, b>0$, $\IG(z\mid a, b)$ represents an inverse-gamma distribution with density $f(z) = b^a/\Gamma(a) z^{-a-1}\exp(-b/z), ~ z>0$.
For $\lambda>0,k>0$, $\Wbll(\lambda,k)$ denotes a Weibull distribution with density function $f(z) = (k/\lambda)(z/\lambda)^{k-1}\exp\{-(z/\lambda)^{k}\}, ~z>0$.
The symbols $\precsim$ and $\succsim$ are used to denote inequalities up to a constant multiple.
The symbol $a \simeq b$ denotes $a \precsim b$ and $a \succsim b$.
For two real numbers $a$ and $b$, $a\vee b$ and $a\wedge b$ denote the bigger and the smaller of the two, respectively.

\newpage
\section{Main Results}\label{sec: main thms}

We consider a DPMM prior $\Pi$ on $f_{X}$, defined as follows.
Set $f_{X} (x) = \sum_{k=1}^{\infty} \pi_{k} ~ \Normal(x\mid \mu_k,h^{2} ) = \int \Normal (x-z\mid 0,h^{2}) dF_{\alpha}(z)  = \{dF_{\alpha}\star \Normal(0,h^{2} )\}(x) $, where $dF_{\alpha}(z) = \sum_{k=1}^{\infty}\pi_{k}\delta_{\mu_k}$ with $\pi_{k} = S_{k} \prod_{\ell=1}^{k-1}(1-S_{\ell})$, $S_{\ell} \sim \Beta(1,\alpha)$, $\mu_{k} \sim P_{0\mu}, h \sim P_{0h}$.
The induced model for the density of $W$ then becomes 
$f_{W} (w) = (f_{X} \star \psi_{\sigma}) (w) = \int \psi_{\sigma} (w-x) f_{X}(x) dx = \sum_{k=1}^{\infty} \pi_{k} \int \Normal(x\mid \mu_k,h^{2})  \psi_{\sigma} (w-x)  dx$.  
We restrict our attention to priors for which $P_{0\mu}$ and $P_{0h}$ satisfy the following conditions.

\begin{Cond} \label{cond: base measure}
1. $P_{0\mu}$ has a positive density on the whole $\rR$. 
2. There exist positive constants $K_{1}$ and $K_{2}$ and a $y$ sufficiently small such that the density of $P_{0h}$, say $f_{0h}$, satisfies $f_{0h}(h) \succsim \exp(-K_{1}h^{-K_{2}})$ for all $h\in (0,y)$.  
\end{Cond}
$\Normal(0,\sigma_{0}^{2})$ satisfies Condition \ref{cond: base measure}.1 and leads to easy posterior computation. 
Condition \ref{cond: base measure}.2 is satisfied by a $\Wbll(\lambda,k)$ or an $ \IG(a,b)$ prior on any positive power of $h$.
A conjugate $\IG(a,b)$ prior on $h^{2}$ may be preferred as it leads to easy posterior computation.

A formal proof of the main consistency results also requires that $f_{0X}$ and $\psi_{\sigma}$ satisfy some minimal assumptions.

\begin{Thm} \label{thm: main thm supersmooth}
Let $2 \leq p \leq \infty$, $f_{0X} \in \C_{p}^{\eta}(\rR)$ and $\Pi$ be the DPMM prior described above. 
Assume there exist constants $c_{1}>0, c_{2} > 1$ and $T > 0$ such that 
\be \label{eq: polynomial tails}
f(z) \leq c_{1} \abs{z}^{-c_{2}} ~~\text{for}~~ \abs{z} \geq T,~~~f \in \{f_{0X}, \psi_{\sigma}\},
\ee 
that is, $f_{0X}$ and $\psi_{\sigma}$ both have polynomially decaying tails.
Also assume $\int \abs{f'} <\infty$ for at least one $f$ in $\{f_{0X}, \psi_{\sigma}\}$. 
Then there exists an $M \in (0,\infty)$ such that for $2 \leq p < \infty$
\be \label{eq: main thm supersmooth}
E_{f_{0W}}\Pi\{ \norm{f_X - f_{0X}}_p \geq M (\log n)^{-\eta/\beta} \mid \bW_{1:n}\} \to 0.
\ee
If further $\int \abs{\phi'_{\sigma}}^{2} \simeq \int \abs{x\psi_{\sigma}(x)}^{2}dx <\infty$,
then (\ref{eq: main thm supersmooth}) also holds for $p = \infty$.
\end{Thm}

We prove Theorem \ref{thm: main thm supersmooth} by verifying a set of sufficient conditions presented in Section \ref{sec: suff thm}.
In density estimation problems consistency results require that the prior assigns sufficient mass in KL balls around $f_{0X}$.
To establish this prior probability bound, it is usually assumed that $f_{0X}$ has exponentially decaying tails. 
See for instance Ghosal and van der Vaart (2007), Kruijer, et al. (2010) and Shen, et al. (2013).
In our density deconvolution problem also we need to establish similar prior concentration bounds, but only for KL balls around $f_{0W}$, the convoluted version of $f_{0X}$. 
And a much weaker polynomial tails condition on $f_{0X}$ and $\psi_{\sigma}$ suffices.
See Section \ref{sec: KL} for details.

\section{Sufficient Conditions for Posterior Convergence} \label{sec: suff thm}
To get posterior contraction rates in general deconvolution problems, we provide a set of sufficient conditions stated in the following theorem.
The proof is based on the proofs of Theorem 2.1 in Ghosal, et al. (2000) and Theorem 3.1 in Ray (2012) and is moved to the Appendix.

\begin{Thm} \label{thm: suff}
Let $\bW_{1:n} = \{W_1, \ldots, W_n\}$ be generated from $f_{0W}= f_{0X} \star \psi_{\sigma}$ where $f_{0X}\in \F$. 
Let $\Pi_{n}$ denote a sequence of priors on $\F$.
Suppose there exist sequences $\xi_{n},\epsilon_{n}$ of positive numbers with $\xi_{n}\to0, \epsilon_{n} \to 0, n\epsilon_{n}^2 \to \infty$,  
a sequence of estimators $f_{n, X}$ based on $\bW_{1:n}$, 
a sequence $\calP_n$ of subsets of $\left\{f_{W}: f_{W}=f_{X} \star \psi_{\sigma},  \norm{E_{f_{W}}f_{n, X} - f_{X}}_p < D\xi_n \right\}$ for some $D>0$, 
and positive constants $C, D_{0}$ and $D_{1}$ such that
\be
\norm{E_{f_{0W}} f_{n, X} - f_{0X}}_p < D_{0}\xi_n,  \label{eq: bias of DKE}  \\
\Pi_{n} \left( \calP_{n}^c\right)  \leq \exp\{-(4+C)n \epsilon_{n}^2\}, \label{eq: suff1} \\
\sup_{f_{W}\in\calP_{n}\cup\{f_{0W}\}}\Pr_{f_{W}}\left( \norm{f_{n, X}  -  E_{f_{W}}f_{n, X}}_p \geq D_{1} \xi_n \right) \leq \exp\{-(4+C)n \epsilon_{n}^2\},  	\label{eq: estconc}  \\
\Pi_{n} \left\{f_W:  \int f_{0W}  \log \frac{f_{0W}}{f_{W}}  < \epsilon_{n}^2,  \int f_{0W} \bigg(\log \frac{f_{0W}}{f_{W}} \bigg)^2 < \epsilon_{n}^2  \right\} \geq \exp(-C n \epsilon_{n}^2). \label{eq: suff2}
\ee
Then there exists an $M \in (0,\infty)$ such that 
\be
E_{f_{0W}}\Pi_{n}( \norm{f_X - f_{0X}}_p \geq M \xi_n \mid \bW_{1:n}) \to 0.
\ee
\end{Thm}

\section{Auxiliary Results}\label{sec: aux results}

\subsection{A Sequence of Estimators $f_{n,X}$}\label{sec: DKE}

For the sequence of estimators $f_{n,X}$, that is used to construct a test function for testing (\ref{eq: testing problem}) 
and appears in the statement of Theorem \ref{thm: suff}, 
we use the deconvoluting kernel density estimator (DKE) (Carroll and Hall, 1988, Stefanski and Carroll, 1990). 
In this subsection we briefly review the DKE and establish condition (\ref{eq: bias of DKE}) of Theorem \ref{thm: suff}.

Taking Fourier transform on both sides of (\ref{eq: deconv}), we have $\phi_{0W}(t) = \phi_{0X}(t)\phi_{\sigma}(t)$.
By Fourier inversion, $f_{0X}$ can thus be obtained as
\bse
f_{0X}(x) = \frac{1}{2\pi} \int_{-\infty}^{\infty} \exp(-itx) \frac{\phi_{0W}(t)}{\phi_{\sigma}(t)} dt.
\ese
Here $\phi_{0W}(t)$ can be estimated by $\phi_{n,W}(t) = n^{-1} \sum_{j=1}^n \exp(itW_j)$, the empirical characteristic function. 
However, since $\phi_{n,W}(t)$ is not a good estimate of  $\phi_{0W}(t)$, particularly at high-frequencies,  
the idea is to incorporate a damping factor $\phi_{K}(h_n t)$ , where $\phi_{K}$, with $\phi_{K}(0) =1$, 
is the Fourier transform of a kernel $K$, and $h_n \downarrow 0$ so that the damping factor $\phi_{K}(h_n t) \to 1$. 
Thus, a DKE is obtained as  
\be\label{eq: est}
f_{n, X}(x) = \frac{1}{2\pi} \int_{-\infty}^{\infty} \exp(-itx) \phi_{n,W}(t)  \frac{\phi_{K}(h_n t)}{\phi_{\sigma}(t)} dt.
\ee 
See Fan (1991) and the references therein for more details. 
We will choose an appropriate $h_n$ in the sequel.  
Under assumptions integrability, (\ref{eq: est}) can be written as a kernel type estimate 
\be\label{eq: estker}
f_{n, X}(x)  = \frac{1}{n} \sum_{j=1}^n \frac{1}{h_n}K_n \bigg( \frac{x - W_j}{h_n} \bigg),
\ee
where 
\be\label{eq: ker}
K_n(x) =  \frac{1}{2\pi} \int_{-\infty}^{\infty} \exp(-itx) \frac{\phi_{K}(t)}{\phi_{\sigma}(t/h_n)} dt.
\ee
Here $h_n$ plays the role of a bandwidth in the kernel $K_n$ in (\ref{eq: estker}). 
It is important to note that $h_n$ also appears in the construction of $K_n$ in (\ref{eq: ker}).  
Note that
\be
&&\hspace{-1cm} E_{f_{0W}} f_{n, X} (x) =  E_{f_{0W}} h_n^{-1} K_n \{(x - W_1)/h_n \}  \nonumber \\
&& = \int f_{0W}(w) h_n^{-1} K_n \{(x-w)/h_n\} dw 
= (f_{0W} \star \wt{K}_{n,W}) (x)  \nonumber \\
&& = \int f_{0W}(w) \frac{1}{2\pi h_{n}} \int \exp\{-it(x-w)/h_{n}\} \frac{\phi_{K}(t)}{\phi_{\sigma}(t/h_n)} dt  ~ dw \nonumber\\
&& = \frac{1}{2\pi h_{n}} \int \exp(-itx/h_{n})\phi_{0W}(t/h_{n}) \frac{\phi_{K}(t)}{\phi_{\sigma}(t/h_n)} dt  \nonumber\\
&& = \frac{1}{2\pi h_{n}} \int \exp(-itx/h_{n})\phi_{0X}(t/h_{n}) \phi_{K}(t) dt \nonumber\\
&& = \frac{1}{2\pi} \int \exp(-itx)\phi_{0X}(t) \phi_{K}(th_{n}) dt  
= (f_{0X} \star \wt{K}_{n,X}) (x),
\ee
where $\wt{K}_{n,W} (z) =  h_{n}^{-1} K_{n}(z/h_n)$ and $\wt{K}_{n,X} (z) =  h_{n}^{-1} K(z/h_n)$. 
Let the kernel $K$ be a $\lfloor\eta\rfloor\th$ order kernel that satisfy the following conditions.
\begin{Cond} \label{cond: kernel}
1. $K(z) = K(-z)$. 
2. $\int z^{r} K(z) dz = 0$ for $r=1,2,\ldots,\lfloor\eta\rfloor$.\\
3. $\int \abs{z}^{\lfloor\eta\rfloor+(\eta-\lfloor\eta\rfloor) p}  \abs{K(z)} dz < \infty$.
\end{Cond}
Examples of higher order kernels can be found in Prakasa Rao (1983).
General methods of construction of higher order kernels starting with simple second order kernels can also be found in Fan and Hu (1992).   
Then, for any $f_{0X} \in C_{p}^{\eta}(\rR)$, using Taylor series expansion, we have
\bse
&& \hspace{-1cm} E_{f_{0W}} f_{n,X} (x) - f_{0X}(x)  =    \int \{f_{0X}(x-h_{n}y) - f_{0X}(x) \}  K (y) dy    \nonumber\\
&&= \int \int_{0}^{1} \frac{(1-t)^{\lfloor \eta \rfloor-1}}{(\lfloor \eta \rfloor-1)!}  \left\{f_{0X}^{(\lfloor \eta \rfloor)}(x-th_{n}y)-f_{0X}^{(\lfloor \eta \rfloor)}(x)\right\} (h_{n}y)^{\lfloor \eta \rfloor}  K(y) dt dy    \nonumber\\
\Rightarrow && \abs{E_{f_{0W}} f_{n, X} - f_{0X}} \precsim h_{n}^{\lfloor \eta \rfloor}  \int \int_{0}^{1} (1-t)^{\lfloor \eta \rfloor-1} \abs{f_{0X}^{(\lfloor \eta \rfloor)}(x-th_{n}y)-f_{0X}^{(\lfloor \eta \rfloor)}(x)} \abs{y}^{\lfloor \eta \rfloor}  \abs{K(y)} dt dy.      
\ese
Let $g(t,y) = (1-t)^{\lfloor \eta \rfloor-1}\abs{y}^{\lfloor \eta \rfloor}  \abs{K(y)}$ and $C =\int \int_{0}^{1} g(t,y) dt dy <\infty$.
Then $g^{\star}(t,y) = C^{-1} g(t,y)$ is a probability density.
Therefore, for $2 \leq p < \infty$, applying Jensen's inequality and Fubini's theorem, we have
\be
&& \hspace{-1cm} \norm{E_{f_{0W}} f_{n, X} - f_{0X}}_{p} \precsim h_{n}^{\lfloor \eta \rfloor}   \left[\int \left\{\int  \int_{0}^{1} \abs{f_{0X}^{(\lfloor \eta \rfloor)}(x-th_{n}y)-f_{0X}^{(\lfloor \eta \rfloor)}(x)} g^{\star}(t,y) dt dy \right\} ^{p} dx\right]^{1/p}   \nonumber\\
&& \leq h_{n}^{\lfloor \eta \rfloor}   \left\{\int \int  \int_{0}^{1} \abs{f_{0X}^{(\lfloor \eta \rfloor)}(x-th_{n}y)-f_{0X}^{(\lfloor \eta \rfloor)}(x)}^{p} g^{\star}(t,y) dx dt dy\right\}^{1/p}   \nonumber\\
&& \precsim h_{n}^{\lfloor \eta \rfloor}   \left\{\int  \int_{0}^{1} \abs{th_{n}y}^{(\eta-\lfloor\eta\rfloor)p} g^{\star}(t,y) dt dy\right\}^{1/p}   \nonumber\\
&& \precsim h_{n}^{\eta}  \left\{\int  \int_{0}^{1} (1-t)^{\lfloor\eta\rfloor-1} t^{(\eta-\lfloor\eta\rfloor)p} \abs{y}^{\lfloor\eta\rfloor+(\eta-\lfloor\eta\rfloor)p} \abs{K(y)} dt dy\right\}^{1/p}   \nonumber\\
&& \precsim h_{n}^{\eta}  \left\{\int \abs{y}^{\lfloor\eta\rfloor+(\eta-\lfloor\eta\rfloor)p} \abs{K(y)} dy\right\}^{1/p} = O(h_{n}^{\eta}).
\ee
Similarly, we also have
\be
&& \norm{E_{f_{0W}} f_{n, X} - f_{0X}}_{\infty} =  \sup_{x} \abs{E_{f_{0W}} f_{n,X} (x) - f_{0X}(x)} \nonumber\\
&& \leq \sup_{x} \int \abs{f_{0X}(x-y) -f_{0X}(x)}h_{n}^{-1}\abs{K(y/h_{n})} dy \leq \frac{C}{\lfloor\eta\rfloor!}\int \abs{y}^{\eta} h_{n}^{-1}\abs{K(y/h_{n})} dy \nonumber\\
&& = \frac{Ch_{n}^{\eta}}{\lfloor\eta\rfloor!}\int \abs{z}^{\eta} \abs{K(z)} dz = O(h_{n}^{\eta}).
\ee

\subsection{Concentration Results for $f_{n,X}$} \label{sec: concentration bounds for DKE}

A main step in the any posterior convergence theorem is the construction  of nonparametric tests for $L^{p}$-alternatives, $2 \leq p \leq \infty$, 
that have sufficiently good exponential bounds on the type I and II errors.  
The degree of concentration of around its expectation in $L_{p}$ norm depends on $p$ 
and has been worked out in Gin\'{e} and Nickl (2012) for usual linear kernel type estimators and wavelets. 
These results have been derived using the Talagrand's inequality (Talagrand, 1996) for empirical processes.  
Concentration bounds for $2 \leq p < \infty$ are obtained in subsection \ref{sec: concentration bounds Talagrand} by adapting these results for deconvolution estimators.  
Concentration bounds for the case $p= \infty$ are obtained separately in subsection \ref{sec: concentration bound for p = infty} using a result due to Dvoretzsky, et al. (1956).
In what follows, the sieve $\calP_{n}$ is the set of all prior realizations convoluted with $\psi_{\sigma}$.

\subsubsection{Concentration Bounds for $2 \leq p < \infty$}  \label{sec: concentration bounds Talagrand}

Let $Z_{1}, \ldots, Z_{n}$ be independently and identically distributed with law $P$  on a measurable space $(S, \S)$. 
Let $\wh{F}_{n,Z} = n^{-1}\sum_{j=1}^{n}\delta_{Z_{i}}$. 
Also let $\G$ be a $P$-centered (i.e., $\int g dP = 0$ for all $g \in \G$) countable class of real-valued measurable functions in $S$ 
and set $\norm{n\wh{F}_{n,Z}}_{\G} = \sup_{g \in \G} \abs{\int g~n~d\wh{F}_{n,Z}}  = \sup_{g \in \G} \abs{\sum_{j=1}^{n}g(Z_{j})}$.  
Let $K_{1}$ and $K_{3}$ be constants such that $\G$ is uniformly bounded by $K_{1}$ and $\sup_{g \in \G}  E g^2(Z) \leq K_{3}$ and set $K_{2}= nK_{3} + 2K_{1} E\norm{n\wh{F}_{n,Z}}_{\G}$. 
Then, Bousquet's (Bousquet, 2003) version of Talagrand's inequality, with constants, is as follows: for every $z \geq 0, n \in \mathbb{N}$, 
\be \label{eq: tallagrand}
\textstyle P \bigg\{ \norm{n\wh{F}_{n,Z}}_{\G}  \geq E  \norm{n\wh{F}_{n,Z}}_{\G} + (2K_{2}z)^{1/2} + K_{1}z/3 \bigg\} \leq 2 e^{-z}.
\ee
This applies to our situation as follows. 
Consider  
\bse
n \{f_{n, X}(x) - E_{f_{W}} f_{n,X}(x)\}  = \sum_{j=1}^n \{\wt{K}_n (x - W_j) -   E_{f_{W}} \wt{K}_{n} (x- W_j)\}
\ese
where  $\{\wt{K}_{n} (x-W_{j}) - E_{f_{W}}\wt{K}_{n} (x-W_{j}), j=1, \ldots, n \}$ are $L_{p}(\rR)$ valued random variables. 
Let $2 \leq p < \infty$ and $q$ be conjugate to $p$ in the sense that $1/p + 1/q = 1$.  
Then, by Hahn-Banach theorem, the separability of $L_{p}$ implies that 
there exists a countable subset $B_{0q}$ of the unit ball $L_{q}$ such that for any $s\in L_{p}$
\be \label{eq:hb}
\norm{s}_{p} =  \sup_{m \in B_{0q}}  \left \vert \int_{\rR} s(x) m(x) dx \right\vert .
\ee
By (\ref{eq:hb}) it is easy to see that 
\be
&& \hspace{-1cm} n \norm{f_{n,X} - E_{f_{W}} f_{n,X}}_{p} = n \sup_{m\in B_{0q}} \left\vert \int \{f_{n,X}(x) - E_{f_{W}} f_{n,X}(x)\} m(x) dx \right\vert   \nonumber\\
&& = \sup_{m\in B_{0q}} \left\vert \int \sum_{j=1}^n \{\wt{K}_n (x - W_j) - E_{f_{W}} \wt{K}_{n} (x- W_j)\}  m(x)dx \right\vert  \nonumber\\
&& = \sup_{m \in B_{0q}} \left\vert \sum_{j=1}^n \int \{\wt{K}_n (x - W_j) - E_{f_{W}} \wt{K}_{n} (x- W_j)\} m(x)dx \right\vert  \nonumber\\
&& = \sup_{k \in \K_{p}} \left\vert \sum_{j=1}^n k(W_{j})\right\vert  = \norm{n\wh{F}_{n,W}}_{\K_{p}},
\ee
where $\K_{p} = \{z \mapsto \int \wt{K}_{n} (x-z) m(x) dx - \int E_{f_{W}} \{ \wt{K}_{n} (x-W_{1})\} m(x) dx :  m \in B_{0q}\}$ and $\wh{F}_{n,W} = n^{-1} \sum_{j=1}^n \delta_{W_{j}}$. 
Note that the class $\K_{p}$ equals to $\K_{p} = \{t  \mapsto (\wt{K}_{n} \star m)(t)  -  E_{f_{W}} \{ (\wt{K}_{n} \star m)(W_{1})\}:  m \in B_{0q}\}$.


To apply (\ref{eq: tallagrand}), we need to find suitable bounds for the envelop 
$K_{1} \geq \sup_{k \in \K_{p}} \abs{k}$ and weak variances $K_{2} \geq \sup_{k \in \K_{p}} E_{f_{W}} k^2(W)$. 
The following three lemmas provide the required bounds.
To establish the bounds, we assume, as in Fan (1991b), that the kernel $K$ in the DKE additionally satisfies the following condition.

\begin{Cond}
$\phi_{K}(t) =0$ for $\abs{t}>1$.
\end{Cond}

\noindent 
For $2 \leq p < \infty$ define $\varpi_{p}(\beta_{0}) = -(1-1/p)$ if $\beta_{0} \geq 0$ and $\varpi_{p}(\beta_{0}) = \{\beta_{0}-(1-1/p)\}$ if $\beta_{0} < 0$.
Also define $\varpi_{\infty}(\beta_{0}) = -1$ if $\beta_{0} \geq 0$ and $\varpi_{\infty}(\beta_{0}) = (\beta_{0}-1)$ if $\beta_{0} < 0$.
Note that for $2 \leq p < \infty$, we have $\norm{\wt{K}_{n}}_{p}^{p}  = h_{n}^{-(p-1)} \norm{K_{n}}_{p}^{p}$ and 
$\norm{\wt{K}_{n}^{2}}_{p}^{p}  = h_{n}^{-(2p-1)} \norm{K_{n}^{2}}_{p}^{p} = h_{n}^{-(2p-1)} \norm{K_{n}}_{2p}^{2p}$.

\begin{Lem} \label{lem: vb1}
$\sup_{k \in \K_{p}} \abs{k}  \precsim K_{1} = h_{n}^{\varpi_{\infty}(\beta_{0})} \exp(h_{n}^{-\beta}/ \varrho)$ for $2 \leq p < \infty$.
\end{Lem}
\begin{Proof}
By Young's inequality, we obtain 
\bse
\norm{(\wt{K}_{n} \star m) (t)}_{\infty} \leq \norm{\wt{K}_{n}}_{p} \norm{m}_{q}.
\ese
Next, we find an upper bound to $\norm{\wt{K}_{n}}_{p}$. 
For $2 \leq p < \infty$, we have $q \in (1,2]$.
Applying Hausdorff-Young's inequality, we have 
\bse
\norm{\wt{K}_{n}}_{p}^{p}  &=&  h_{n}^{-(p-1)} \norm{K_n}_{p}^{p}  = h_{n}^{-(p-1)}  \int \left\vert \frac{1}{2\pi} \int\exp(-itz) \frac{\phi_{K}(t)}{\phi_{\sigma}(t/h_n)} dt \right\vert^{p}dz 
\precsim h_{n}^{-(p-1)} \norm{ \frac{\phi_{K}(\cdot)}{\phi_{\sigma}(\cdot/h_{n})}}_{q}^{p} \\
\Rightarrow \norm{\wt{K}_{n}}_{p}  &=&  h_{n}^{-(1-1/p)} \norm{K_n}_{p} \precsim h_{n}^{-(1-1/p)} \norm{ \frac{\phi_{K}(\cdot)}{\phi_{\sigma}(\cdot/h_{n})}}_{q}.
\ese
Note that
\bse
&& \norm{ \frac{\phi_{K}(\cdot)}{\phi_{\sigma}(\cdot/h_{n})}}_{q}^{q} =  \int_{\abs{t} \leq  M h_{n}} \abs{ \frac{\phi_{K}(t )}{\phi_{\sigma}(t/h_{n})}}^{q} dt  +  \int_{M h_{n} <\abs{t} \leq 1} \abs{ \frac{\phi_{K}(t)}{\phi_{\sigma}(t/h_{n})}}^{q} dt \\
&\leq&  h_{n} \int_{\abs{z} \leq M} \abs{ \frac{\phi_{K}(z h_n)}{\phi_{\sigma}(z)}}^{q} dz  +  (2/d_{0})^{q} h_{n}^{q\beta_{0}}\int_{M h_{n} <\abs{t} \leq 1}  \exp(q \abs{t}^{\beta} h_{n}^{-\beta}/ \varrho)  \abs{t}^{-q\beta_{0}} \abs{\phi_{K}(t)}^{q}  dt \\
&\leq & \frac{\sup \abs{\phi_{K}}^{q}}{\inf_{z \leq M}\abs{\phi_{\sigma}(z)}^{q}}   M h_{n}  +   (2/d_{0})^{q}h_{n}^{q\beta_{0}} \exp(q h_{n}^{-\beta}/ \varrho) \int_{M h_{n} <\abs{t} \leq 1} \abs{t}^{-q\beta_{0}} \abs{\phi_{K}(t)}^{q}  dt.
\ese
When $\beta_{0}\geq 0$, using $\abs{t} \geq M h_{n} \Leftrightarrow \abs{t}^{-q\beta_{0}} \leq (Mh_{n})^{-q\beta_{0}}$ for the second term, we have
\bse
 \norm{ \frac{\phi_{K}(\cdot)}{\phi_{\sigma}(\cdot/h_{n})}}_{q}^{q}    &\precsim&   \frac{\sup \abs{\phi_{K}}^{q}}{\inf_{z \leq M}\abs{\phi_{\sigma}(z)}^{q}}  M h_{n}  +   (2/d_{0})^{q}M^{-q\beta_{0}} \exp(q h_{n}^{-\beta}/ \varrho) \int  \abs{\phi_{K}(t)}^{q}  dt  \\
&\precsim& \exp(q h_{n}^{-\beta}/ \varrho)  \\
\Rightarrow \norm{ \frac{\phi_{K}(\cdot)}{\phi_{\sigma}(\cdot/h_{n})}}_{q}  &\precsim& \exp(h_{n}^{-\beta}/ \varrho)
\ese
When $\beta_{0}<0$, using $ \abs{t}^{-q\beta_{0}} \leq 1$, we have
\bse
\norm{ \frac{\phi_{K}(\cdot)}{\phi_{\sigma}(\cdot/h_{n})}}_{q}^{q}  &\precsim&   \frac{\sup \abs{\phi_{K}}^{q}}{\inf_{z \leq M}\abs{\phi_{\sigma}(z)}^{q}}   M h_{n}  +   (2/d_{0})^{q} h_{n}^{q\beta_{0}} \exp(q h_{n}^{-\beta}/ \varrho) \int \abs{\phi_{K}(t)}^{q}  dt  \\
&\precsim& h_{n}^{q\beta_{0}} \exp(q h_{n}^{-\beta}/ \varrho)  \\
\Rightarrow \norm{ \frac{\phi_{K}(\cdot)}{\phi_{\sigma}(\cdot/h_{n})}}_{q} &\precsim& h_{n}^{\beta_{0}} \exp(h_{n}^{-\beta}/ \varrho)  \\
\ese
Combining, we have
\bse
\norm{\wt{K}_{n}}_{p}  &\precsim&  h_{n}^{\varpi_{p}(\beta_{0})} \exp(h_{n}^{-\beta}/ \varrho)
\ese
For {\bf $p=\infty$}, proceeding in a similar manner, we have
\bse
\norm{\wt{K}_{n}}_{\infty}  &\precsim&   h_{n}^{\varpi_{\infty}(\beta_{0})} \exp(h_{n}^{-\beta}/ \varrho)
\ese
A general bound for $2 \leq p < \infty$ is thus obtained for $p=\infty$ and Lemma  \ref{lem: vb1} follows.
\end{Proof}

\begin{Lem}  \label{lem: vb2}
 $\sup_{k \in \K_{p}} E_{f_{W}} k^2(W) \precsim K_{3} = h_{n}^{2\varpi_{\infty}(\beta_{0})} \exp(2 h_{n}^{-\beta}/ \varrho)$ for $2 \leq p < \infty$.
\end{Lem}
\begin{Proof}
By H\"{o}lder's inequality and Young's inequality, 
\bse
 E_{f_{W}} \{ (\wt{K}_n \star m) (W)\}^2  &=&  \int \{ (\wt{K}_n \star m) (w)\}^2 f_{W}(w)dw \leq \norm{f_{W}}_{1} \norm{(\wt{K}_{n} \star m)^2}_{\infty}  \\
&=& \norm{\wt{K}_{n} \star m }_{\infty}^{2}  \leq  \norm{\wt{K}_{n}}_{p}^{2} \norm{m}_{q}^{2} .
\ese
As in Lemma \ref{lem: vb1}, a common bound is obtained for $p=\infty$ and Lemma \ref{lem: vb2} follows.
\end{Proof}

\begin{Lem}\label{lem: vb3}
$E_{f_{W}}\norm{n \wh{F}_{n,W}}_{\K_{p}} \precsim  n^{1/2} h_{n}^{\varpi_{p}(\beta_{0})} \exp(h_{n}^{-\beta}/\varrho)$ for $2 \leq p < \infty$.  
\end{Lem}

\begin{Proof}
Let $ \{\wt{K}_n (x - w) - E_{f_{W}} \wt{K}_{n} (x- W_{1})\} = H(x,w)$.
Then
\bse
 && \hspace{-1cm} E_{f_{W}}\norm{n \wh{F}_{n,W}}_{\K_{p}} = E_{f_{W}} n \norm{f_{n,X} - E_{f_{W}} f_{n,X}}_{p}  = E_{f_{W}} \norm{\sum_{j=1}^{n}H(\cdot,W_{j})}_{p}.
\ese

\noindent {\bf Case $2 < p<\infty$:}
Using Jensen's inequality we have
\bse
&&\hspace{-1cm} E_{f_{W}}\norm{\sum_{j=1}^{n}H(\cdot,W_{j})}_{p} = E_{f_{W}}\left\{\int \abs{\sum_{j=1}^{n}H(x,W_{j})}^{p} dx\right\}^{1/p}  \\
&\leq& \left\{E_{f_{W}} \int \abs{\sum_{j=1}^{n}H(x,W_{j})}^{p} dx\right\}^{1/p}  
=  \left\{\int E_{f_{W}} \abs{\sum_{j=1}^{n}H(x,W_{j})}^{p} dx\right\}^{1/p}. 
\ese
By Hoffmann-J{\o}rgensen's inequality we have
\bse
&&\hspace{-1cm}  E_{f_{W}} \abs{\sum_{j=1}^{n}H(x,W_{j})}^{p}    \leq   C_{p}^{p} \left[E \abs{\sum_{j=1}^{n}H(x,W_{j})} + \left\{E_{f_{W}} \max_{1\leq j \leq n} \abs{H(x,W_{j})}^{p}\right\}^{1/p} \right]^{p}   \\
&& \leq C_{p}^{p} 2^{p-1} \left[\left\{E \abs{\sum_{j=1}^{n}H(x,W_{j})}\right\}^{p} + E_{f_{W}} \max_{1\leq j \leq n} \abs{H(x,W_{j})}^{p} \right], 
\ese
where $C_{p} = kp/(1+\log p)$. 
\bse
&&\hspace{-1cm}  E_{f_{W}} \abs{\sum_{j=1}^{n}H(x,W_{j})}  \leq \left[ E_{f_{W}} \left\{\sum_{j=1}^{n}H(x,W_{j})\right\}^{2}\right]^{1/2}   =  \left[ \var_{f_{W}} \left\{\sum_{j=1}^{n}H(x,W_{j})\right\}\right]^{1/2}   \\
&& = n^{1/2} \left[\var_{f_{W}} \left\{H(x,W_{1})\right\} \right]^{1/2} = n^{1/2} \left[\var_{f_{W}} \left\{\wt{K}_{n}(x-W_{1})\right\} \right]^{1/2}  \\
&& \leq n^{1/2} \left\{E_{f_{W}} {\wt{K}_{n}^{2}(x-W_{1})} \right\}^{1/2}.
\ese
Also
\bse
E_{f_{W}} \max_{1\leq j \leq n} \abs{H(x,W_{j})}^{p} \leq 2^{p-1} \left[ E_{f_{W}} \max_{1\leq j \leq n} \abs{\wt{K}_{n}(x-W_{j})}^{p} + \left\{E_{f_{W}} \abs{\wt{K}_{n}(x-W_{j})}\right\}^{p}\right].
\ese
Combining we have
\bse
&&\hspace{-1cm}  E_{f_{W}} \abs{\sum_{j=1}^{n}H(x,W_{j})}^{p}   \leq   C_{p}^{p} 2^{p-1} \left[\left\{E \abs{\sum_{j=1}^{n}H(x,W_{j})}\right\}^{p} + E_{f_{W}} \max_{1\leq j \leq n} \abs{H(x,W_{j})}^{p} \right]    \nonumber\\
&&  \leq   C_{p}^{p} 2^{p-1} \left[(n^{p/2}+2^{p-1}) \left\{E_{f_{W}} \abs{\wt{K}_{n}^{2}(x-W_{1})} \right\}^{p/2}+ 2^{p-1} E_{f_{W}} \max_{1\leq j \leq n} \abs{\wt{K}_{n}(x-W_{j})}^{p} \right].
\ese
Therefore,
\be
&&\hspace{-1cm} E_{f_{W}}\norm{\sum_{j=1}^{n}H(\cdot,W_{j})}_{p}  \precsim \left[ \int n^{p/2} \left\{E_{f_{W}}\abs{\wt{K}_{n}^{2} (x - W_{1})} \right\}^{p/2} dx  + \int E_{f_{W}} \max_{1\leq j \leq n} \abs{\wt{K}_{n} (x - W_{j})}^{p}  dx \right]^{1/p}  \nonumber\\
&&\leq \left[ n^{p/2}\int \left(\abs{\wt{K}_{n}^{2}} \star f_{W}\right)^{p/2}(x) dx  + \int E_{f_{W}} \sum_{j=1}^{n} \abs{\wt{K}_{n} (x - W_{j})}^{p}  dx \right]^{1/p}  \nonumber\\
&&= \left[ n^{p/2}\int \left(\abs{\wt{K}_{n}^{2}} \star f_{W}\right)^{p/2}(x) dx  + n \int E_{f_{W}} \abs{\wt{K}_{n} (x - W_{j})}^{p}  dx \right]^{1/p}  .\label{eq: Talagrand expectation norm 1}
\ee
The first term in (\ref{eq: Talagrand expectation norm 1}) can be bounded by Young's inequality and Hausdorff-Young's inequality
as follows.
\bse
&& \hspace{-1cm}  \left\{\int \left(\abs{\wt{K}_{n}^{2}} \star f_{W}\right)^{p/2}(x) dx \right\}^{2/p}  = \norm{\abs{\wt{K}_{n}^{2}} \star f_{W}}_{p/2}  \leq  \norm{\wt{K}_{n}^{2}}_{p/2} \norm{f_{W}}_{1} = h_{n}^{-(2-2/p)} \norm{K_{n}^{2}}_{p/2} \\
&& = h_{n}^{-(2-2/p)}  \norm{K_{n}}_{p}^{2} \precsim h_{n}^{-(2-2/p)}   \norm{\frac{\phi_{K}(\cdot)}{\phi_{\sigma}(\cdot /h_n)}}_{q}^{2},
\ese
where $q = p/(p-1) \in (1, 2]$.
As in the proof of Lemma \ref{lem: vb1}, with $\varpi' = 0$ if $\beta_{0} \geq 0$ and $\varpi' = \beta_{0}$ if $\beta_{0} < 0$, 
we have $\norm{\phi_{K}(\cdot)/\phi_{\sigma}(\cdot /h_n)}_{q} \precsim h_{n}^{\varpi'} \exp(h_{n}^{-\beta}/\varrho)$. 
Note that $\{p\varpi'-(p-1)\} = p\varpi_{p}(\beta_{0})$.
Therefore, 
\bse
&& \hspace{-1cm}  \left\{\int \left(\abs{\wt{K}_{n}^{2}} \star f_{W}\right)^{p/2}(x) dx \right\}^{2/p}  
\precsim h_{n}^{-(2-2/p)}   \norm{\frac{\phi_{K}(\cdot)}{\phi_{\sigma}(\cdot /h_n)}}_{q}^{2}   \precsim h_{n}^{2\varpi'-(2-2/p)}  \exp(2h_{n}^{-\beta}/\varrho) \\
\Rightarrow && \int \left(\abs{\wt{K}_{n}^{2}} \star f_{W}\right)^{p/2}(x) dx  \precsim h_{n}^{p\varpi'-(p-1)}  \exp(ph_{n}^{-\beta}/\varrho) = h_{n}^{p\varpi_{p}(\beta_{0})}  \exp(ph_{n}^{-\beta}/\varrho).
\ese
Using Young's inequality and Hausdorff-Young inequality, the second term in (\ref{eq: Talagrand expectation norm 1}) can be bounded as follows.
\bse
&& \hspace{-1cm} \int E_{f_{W}} \abs{\wt{K}_{n} (x - W_{j})}^{p}  dx  = \int \int\abs{\wt{K}_{n} (x -w)}^{p} f_{W}(w)~dw ~dx  \\
&& = \int \left(\abs{\wt{K}_{n}}^{p} \star f_{W} \right) (x) ~dx  = \norm{\abs{\wt{K}_{n}}^{p} \star f_{W}}_{1}  \leq \norm{\abs{\wt{K}_{n}}^{p}}_{1} \norm{f_{W}}_{1} \\
&& =  \norm{\wt{K}_{n}}_{p}^{p}  = h_{n}^{-(p-1)}\norm{K_{n}}_{p}^{p} \precsim h_{n}^{-(p-1)} \norm{\frac{\phi_{K}(\cdot)}{\phi_{\sigma}(\cdot /h_n)}}_{q}^{p}   \\
&& \precsim h_{n}^{p\varpi'-(p-1)} \exp(ph_{n}^{-\beta}/\varrho) = h_{n}^{p\varpi_{p}(\beta_{0})}  \exp(ph_{n}^{-\beta}/\varrho). 
\ese
Therefore,
\bse
&&\hspace{-1cm}  n E_{f_{W}} \norm{f_{n,X} - E_{f_{W}} f_{n,X}}_{p}  
\precsim \left[ n^{p/2} h_{n}^{p\varpi_{p}(\beta_{0})} \exp(ph_{n}^{-\beta}/\varrho) + n h_{n}^{p\varpi_{p}(\beta_{0})} \exp(ph_{n}^{-\beta}/\varrho) \right]^{1/p}  \nonumber\\
&& \precsim  n^{1/2} h_{n}^{\varpi_{p}(\beta_{0})} \exp(h_{n}^{-\beta}/\varrho).
\ese
The {\bf case p = 2} is similar but much simpler. Hence part 2 of Lemma \ref{lem: vb3} follows.

\end{Proof}

\subsubsection{Concentration Bound for $p=\infty$} \label{sec: concentration bound for p = infty}
Next we establish concentration bounds for the DKE in $L_{\infty}$ norm.
We assume that the kernel $K$ satisfies the following additional conditions.

\begin{Cond}
1. $\sup_{\abs{t}\leq M} \abs{\phi'_{K}(t)}^{2} < \infty$,
2. $\int \abs{\phi'_{K}}^{2} < \infty$.
\end{Cond}

According to Dvoretzsky, et al. (1956) there exists positive constants $C_{1}$ and $C_{2}$, $C_{2}\in(0,2]$, such that for any $\lambda>0$ 
\be
\Pr\left(\norm{F_{W} - F_{n,W}}_{\infty} \geq \lambda n^{-1/2}\right)  \leq C_{1}\exp(-\alpha \lambda^2).
\ee
Let $V_{n} = V(K_{n}) = \int \abs{K'_{n}(x)}dx$, the total variation of $K_n$.
Integration by parts gives
\bse
&& \hspace{-1cm} \norm{E_{f_{W}} f_{n, X} - f_{n,X}}_{\infty} =  \sup_{x} \left\vert \frac{1}{h_n} \int K_{n} \left(\frac{x-y}{h_n}\right) dF_{W}(y)   - \frac{1}{h_n} \int K_{n} \left(\frac{x-y}{h_n}\right) dF_{n,W}(y)  \right\vert    \nonumber\\
&&\leq \sup_{x} \left\{ \frac{1}{h_{n}}  \int \abs{F_{W}(y)-F_{n,W}(y)} \abs{\frac{1}{h_n}K'_{n} \left(\frac{x-y}{h_n}\right)} dy \right\}  \nonumber\\
&&\leq  \frac{1}{h_{n}}   \sup_{x} \abs{F_{W}(x) - F_{n,W}(x)}  \int \abs{\frac{1}{h_n}K'_{n} \left(\frac{x-y}{h_n}\right)dy }  \nonumber\\
&& \leq \frac{V_n}{h_n}  \sup_{x} \abs{F_{W}(x) - F_{n,W}(x)}.
\ese
Therefore, we have
\be
&& \hspace{-2cm}  \Pr\left(\norm{E_{f_{W}}f_{n,X} - f_{n,X}}_{\infty} \geq D_1 \xi_{n} \right) \leq \Pr\left(\norm{F_{W} - F_{n,W}}_{\infty} \geq  h_{n}D_{1}\xi_{n}/V_{n} \right)  \nonumber\\
&& \leq C_{1}\exp(-C_{2} nh_{n}^{2}D_{1}^{2}\xi_{n}^{2}/V_{n}^{2}). \label{eq: Dvoretzsky}
\ee
Let $w(x) = 1/(1+x^2)$. Applying Jensen's inequality, we have
\bse
V_{n}^{2} &=& \left\{\int \abs{K'_{n}(x)}dx \right\}^{2} = \left\{\int \abs{K'_{n}(x)}w^{-1}(x)w(x)dx \right\}^{2}  \\
&\precsim& \int \abs{K'_{n}(x)}^{2}w^{-2}(x)w(x)dx  \leq \int (1+x^2) \abs{K'_{n}(x)}^{2}dx.
\ese
It is easy to check that $K'_{n}(x) = \phi_{g_{K}}(x)$, where $g_{K}(t) =   (2\pi)^{-1} it \phi_{K}(t)/\phi_{\sigma}(t/h_{n})$.
Some simple manipulation also gives us $(-it) \phi_{g_{K}}(t) = \phi_{g'_{K}}(t)$.
Application of Parseval's identity then gives
\bse
\int \abs{K'_{n}(x)}^{2} dx &\simeq&  \int \abs{g_{K}(t)}^{2}dt,  \\
\int \abs{x K'_{n}(x)}^{2} dx &=& \int \abs{x \phi_{g_{K}}(x)}^{2} dx \simeq \int \abs{\phi_{g'_{K}}(t)}^{2}dt \simeq \int \abs{g'_{K}(t)}^{2}dt.
\ese
Therefore,
\bse
&& \hspace{-1cm} h_{n}^{-2}V_{n}^{2} \precsim h_{n}^{-2} \int (1+x^2) \abs{K'_{n}(x)}^{2} dx \simeq h_{n}^{-2} \left\{ \int \abs{K'_{n}(x)}^{2} dx +  \int x^2 \abs{K'_{n}(x)}^{2} dx  \right\}   \\
&\precsim& h_{n}^{-2} \left\{\int \abs{g_{K}(t)}^{2}dt + \int \abs{g'_{K}(t)}^{2} dt\right\}  \\
&\precsim& h_{n}^{-2} \left\{\int \abs{\frac{t\phi_{K}(t)}{\phi_{\sigma}(t/h_{n})}}^{2}dt + \int \abs{\frac{t\phi'_{K}(t)}{\phi_{\sigma}(t/h_{n})}}^{2} dt  +  \int \abs{\frac{\phi_{K}(t)}{\phi_{\sigma}(t/h_{n})}}^{2}dt   +  \int \abs{\frac{(t/h_{n})\phi_{K}(t)\phi'_{\sigma}(t/h_{n})}{\phi_{\sigma}^{2}(t/h_{n})}}^{2} dt \right\}
\ese
The first term can be bounded as follows.
\bse
&&\hspace{-1cm} h_{n}^{-2} \int \abs{\frac{t\phi_{K}(t)}{\phi_{\sigma}(t/h_{n})}}^{2}dt = h_{n}^{-2} \int_{\abs{t}\leq Mh_{n}} \abs{\frac{t\phi_{K}(t)}{\phi_{\sigma}(t/h_{n})}}^{2}dt  +  h_{n}^{-2} \int_{Mh_{n}< \abs{t} \leq 1} \abs{\frac{t\phi_{K}(t)}{\phi_{\sigma}(t/h_{n})}}^{2}dt \\
&\precsim&  h_{n}  \sup_{\abs{z}\leq M} \abs{\frac{z}{\phi_{\sigma}(z)}}^{2} \sup\abs{\phi_{K}(z)}^{2} M  +  h_{n}^{2\beta_{0}-2} \exp(2h_{n}^{-\beta}/\varrho) \int_{Mh_{n}< \abs{t} \leq 1} \abs{t}^{2-2\beta_{0}}\abs{\phi_{K}(t)}^{2}dt  \\
&\precsim& 
\left\{ \begin{array}{cc}  
\exp(2h_{n}^{-\beta}/\varrho) & if ~ \beta_{0} \geq 1\\
h_{n}^{2\beta_{0}-2} \exp(2h_{n}^{-\beta}/\varrho) & if ~ \beta_{0} < 1.
\end{array}\right.
\ese
\noindent The second term can be bounded as follows.
\bse
&&\hspace{-1cm} h_{n}^{-2} \int \abs{\frac{t\phi'_{K}(t)}{\phi_{\sigma}(t/h_{n})}}^{2}dt  = h_{n}^{-2} \int_{\abs{t}\leq Mh_{n}} \abs{\frac{t\phi'_{K}(t)}{\phi_{\sigma}(t/h_{n})}}^{2}dt  +  h_{n}^{-2} \int_{Mh_{n}< \abs{t} \leq 1} \abs{\frac{t\phi'_{K}(t)}{\phi_{\sigma}(t/h_{n})}}^{2}dt \\
&& \precsim  h_{n} \sup_{\abs{z}\leq M} \abs{\frac{z}{\phi_{\sigma}(z)}}^{2} \sup_{\abs{z}\leq M} \abs{\phi'_{K}(z)}^{2} M  +  h_{n}^{2\beta_{0}-2} \exp(2h_{n}^{-\beta}/\varrho) \int_{Mh_{n}< \abs{t} \leq 1} \abs{t}^{2-2\beta_{0}}\abs{\phi'_{K}(t)}^{2}dt  \\
&& \precsim  
\left\{
\begin{array}{cc} 
\exp(2h_{n}^{-\beta}/\varrho) & if ~ \beta_{0} \geq 1,\\
h_{n}^{2\beta_{0}-2} \exp(2h_{n}^{-\beta}/\varrho) & if ~ \beta_{0} < 1.
\end{array}
\right. 
\ese
The third term can be bounded in a similar way as follows.
\bse
&&\hspace{-1cm} h_{n}^{-2} \int \abs{\frac{\phi_{K}(t)}{\phi_{\sigma}(t/h_{n})}}^{2}dt  
\precsim  
\left\{
\begin{array}{cc} 
h_{n}^{-2} \exp(2h_{n}^{-\beta}/\varrho) & if ~ \beta_{0} \geq 0,\\
h_{n}^{2\beta_{0}-2} \exp(2h_{n}^{-\beta}/\varrho) & if ~ \beta_{0} \geq 0.
\end{array}
\right. 
\ese
\noindent The fourth term can be bounded as follows.
\bse
&&\hspace{-1cm} h_{n}^{-2}\int \abs{\frac{(t/h_{n})\phi_{K}(t)\phi'_{\sigma}(t/h_{n})}{\phi_{\sigma}^{2}(t/h_{n})}}^{2} dt \\
&&\leq h_{n}^{-1} \sup_{\abs{z}\leq M} \abs{\frac{z\phi'_{\sigma}(z)}{\phi_{\sigma}(z)}}^{2}  \sup \abs{\phi_{K}(z)}^{2} M \\
&&~~~~~ + (2/d_{0})^{2} h_{n}^{4\beta_{0}-4} \int_{Mh_{n} < \abs{t} \leq 1} \abs{t}^{2-4\beta_{0}}\exp(4\abs{t}^{\beta}h_{n}^{-\beta}/\varrho) \abs{\phi_{K}(t)\phi'_{\sigma}(t/h_{n})}^{2} dt  \\
&& \precsim  h_{n}^{4\beta_{0}-4} \exp(4h_{n}^{-\beta}/\varrho) \int_{Mh_{n} < \abs{t} \leq 1}\abs{t}^{2-4\beta_{0}}\abs{\phi_{K}(t)\phi'_{\sigma}(t/h_{n})}^{2} dt \\
&& \precsim  
\left\{
\begin{array}{cc} 
h_{n}^{-1} \exp(4h_{n}^{-\beta}/\varrho) &  if ~ \beta_{0} \geq 1/2,\\
h_{n}^{4\beta_{0}-4} \exp(4h_{n}^{-\beta}/\varrho) &  if ~ \beta_{0} < 1/2.
\end{array}
\right. 
\ese
Define $\varpi^{\infty}(\beta_{0}) = -1$ if $\beta_{0} \geq 1/2$, and $\varpi^{\infty}(\beta_{0}) = (4\beta_{0}-4)$ if $\beta_{0} < 1/2$.
Combining all the terms we have, 
\bse
\hspace{-1cm}   h_{n}^{-2} V_{n}^{2} \precsim h_{n}^{-2}\int (1+x^2) \abs{K'_{n}(x)}^{2} dx &\precsim&  h_{n}^{\varpi^{\infty}(\beta_{0})}\exp(4 h_{n}^{-\beta}/\varrho).
\ese

\subsection{Kullback-Leibler Conditions}\label{sec: KL}

In this subsection we show that the KL conditions hold under the polynomial tails assumptions (\ref{eq: polynomial tails}) on $f_{0X}$ and $\psi_{\sigma}$. 
The proof follows the arguments of the proof of Theorem 4 of Shen, et al. (2013), 
but requires new calculations and nontrivial adjustments to adapt it to the deconvolution set-up.

We use Lemma 7 of Ghosal and van der Vaart (2007) stated below.   
\begin{Lem} \label{lem: lemma 7 of GvdW}
There exists a $\lambda_0 \in (0,1)$  such that for any two densities $f_0$ and $f_1$ with $\lambda < \lambda_0$ we have 
\bse
&& \int f_0 \log \frac{f_0}{f_1} \leq d_{H}^2(f_0, f_1) \bigg( 1 + 2 \log \frac{1}{\lambda} \bigg) + 2 \int_{f_1/f_0 \leq \lambda} f_{0} \log \frac{f_1}{f_0} \\
&&\int f_0 \bigg(\log  \frac{f_0}{f_1}\bigg)^2 \leq d_{H}^2(f_0, f_1) \bigg( 12 + 2 \log \frac{1}{\lambda} \bigg)^2 + 8 \int_{f_1/f_0 \leq \lambda}  f_{0}\bigg(\log \frac{f_1}{f_0}\bigg)^2. 
\ese
\end{Lem}

Let $\epsilon_{n} = n^{-\gamma}(\log n)^{t}, ~\tau_{n} = \epsilon_{n}^{2}\{\log(1/\epsilon_{n})\}^{-2}, ~\delta_{n} = \epsilon_{n}^{2 b_{1}}$. 
Choose $b_{1}>1$ such that $\epsilon_{n}^{b_1} \{\log(1/\epsilon_{n})\}^{5/4} \leq \epsilon_{n}$. 
Then $\delta_{n}^{1/2} \leq \delta_{n}^{1/2}\log(1/\delta_{n})^{1/4} \precsim \epsilon_{n}^{b_1}\{\log(1/\epsilon_{n})\}^{5/4-1} \leq \epsilon_{n}\{\log(1/\epsilon_{n})\}^{-1} = \tau_{n}^{1/2}$.
Recall that the polynomial tails condition (\ref{eq: polynomial tails}) on $f_{0X}$ and $\psi_{\sigma}$ implies the existence of positive constants 
$c_{1}, c_{2}$ and $T$ such that for $f\in \{f_{0X},\psi_{\sigma}\}$ $f(z) \leq c_{1}\abs{z}^{-c_{2}}$ whenever $\abs{z} \geq T$.
For any two sequences of probability measures $F'_{1,n}$ (compactly supported) and $F'_{2,n}$ (discrete) on $\rR$, we have
\be \label{eq: H}
&& \hspace{-1cm} d_{H}(f_{W}, f_{0W}) \leq d_{H}(f_{0X} \star \psi_{\sigma}, f_{0X} \star \psi_{\sigma} \star \Normal(0,\tau_{n}^{2})) \nonumber\\
&&~ + ~ d_{H}(f_{0X} \star \psi_{\sigma} \star \Normal(0,\tau_{n}^{2}), dF'_{1,n}\star \psi_{\sigma} \star \Normal(0,\tau_{n}^{2}))  \nonumber\\
&&~~ +~ d_{H}(dF'_{1,n} \star \psi_{\sigma} \star \Normal(0,\tau_{n}^{2}), dF'_{2,n}\star \psi_{\sigma} \star \Normal(0,\tau_{n}^{2}))  \nonumber\\
&&~~~ + ~ d_{H} (dF'_{2,n} \star \psi_{\sigma} \star \Normal(0,\tau_{n}^{2}), dF_{\alpha}\star \psi_{\sigma} \star \Normal(0,\tau_{n}^{2})) \nonumber\\
&&~~~~ + ~ d_{H} (dF_{\alpha}\star \psi_{\sigma} \star \Normal(0,\tau_{n}^{2}), dF_{\alpha}\star \psi_{\sigma} \star \Normal(0,h^{2})).
\ee 
We consider bounds for each term on the right hand side of (\ref{eq: H}) one by one.

\noindent {\bf First Term:} 
Since $\int \abs{f'} <\infty$ for at least one $f$ in $\{f_{0X},\psi_{\sigma}\}$, some simple calculations yield 
\be \label{eq: KL1}
\norm{f_{0X} \star \psi_{\sigma} - f_{0X} \star \psi_{\sigma} \star \Normal(0,\tau_{n}^{2})}_{1} \leq \norm{f - f \star \Normal(0,\tau_{n}^{2})}_1 \precsim  \tau_{n}.
\ee

\noindent {\bf Second Term:} 
Define $A_n = [-a_n,a_n]$, $A_{2,n} = [-a_n/2,a_n/2]$ and let $f\in \{f_{0X},\psi_{\sigma}\}$.
The tail conditions on $f$ implies the existence of a small $\gamma_{1}>0$ such that $\int f^{1-\gamma_{1}} (x)dx < \infty$.
Define $A_{3,n} = \{z: f(z)^{-\gamma_{1}} > \tau_{n}^{-2(1+\gamma_{2})}\} = \{z: f(z)< \tau_{n}^{\gamma_{3}}\}$ where $\gamma_{3} = 2(1+\gamma_{2})/\gamma_{1}$. 
Also define $A_{4,n} = \{z: c_{1}\abs{z}^{-c_{2}}) \leq \tau_{n}^{\gamma_{3}}\} = \{z: \abs{z} \geq c_{1}^{1/c_{2}} \tau_{n}^{-\gamma_{3}/c_{2}}\}$
and $A_{5,n} =  \{z: f(z) \leq c_{1}\abs{z}^{-c_{2}}\}$.
For $n$ sufficiently large $c_{1}^{1/c_{2}} \tau_{n}^{-\gamma_{3}/c_{2}} \geq T$.
Then $A_{4,n} \subseteq \{z: \abs{z} \geq T\}  \subseteq  A_{5,n}$.
Take $a_{0} \geq 2c_{1}^{1/c_{2}}$ and $a_{n} = a_{0} \tau_{n}^{-\gamma_{3}/c_{2}}$.
For $n$ large enough, we have $a_{n}/2 = (a_{0}/2)\tau_{n}^{-\gamma_{3}/c_{2}} \geq c_{1}^{1/c_{2}} \tau_{n}^{-\gamma_{3}/c_{2}}$. 
Hence, $A_{2,n}^{c} \subseteq A_{4,n}  =  A_{4,n} \cap A_{5,n} \subseteq  \{z: f(z) \leq \tau_{n}^{\gamma_{3}}\} = A_{3,n}$.
Therefore, using Markov's inequality, for $f\in \{f_{0X},\psi_{\sigma}\}$ we have 
\be \label{eq: tail both}
\Pr_{f} (A_{2,n}^c) \leq \Pr_{f}(A_{3,n}) = \Pr_{f}\{z: f(z)^{-\gamma_{1}} > \tau_{n}^{-2(1+\gamma_{2})}\} \precsim \tau_{n}^{2+2\gamma_{2}}.
\ee

Define $dF'_{1,n} = f_{0X}1_{A_n}/\int_{A_n}f_{0X}(x)dx$. 
Then, using Young's inequality, we have
\be \label{eq: KL2}
&& \hspace{-1cm} \norm{f_{0X} \star \psi_{\sigma} \star \Normal(0,\tau_{n}^{2}) - dF'_{1,n} \star \psi_{\sigma} \star \Normal(0,\tau_{n}^{2})}_{1} \nonumber\\
&& \leq \norm{\psi_{\sigma} \star \Normal(0,\tau_{n}^{2})}_{1} \norm{f_{0X} - dF'_{1,n}}_1 = \int \abs{f_{0X}(x) - dF'_{1,n}(x)} dx  \nonumber\\
&& =  \int_{A_n} \abs{f_{0X}(x) - dF'_{1,n}(x)} dx +  \int_{A_n^c} \abs{f_{0X}(x) - dF'_{1,n}(x)} dx  \nonumber\\
&& =  \int_{A_n} \{1/\textstyle{\int_{A_n}f_{0X}(z)dz} - 1\} f_{0X}(x) dx +  \int_{A_n^c} f_{0X}(x) dx  \nonumber\\
&& =  1-\Pr_{f_{0X}}(A_n) + \Pr_{f_{0X}}(A_n^c)   \nonumber\\
&& =  2 ~ \Pr_{f_{0X}}(A_n^c) \leq 2 ~ \Pr_{f_{0X}} (A_{2,n}^c) \precsim  \tau_{n}^{2+2\gamma_{2}}.
\ee

\noindent {\bf Third Term:} 
%
Using Lemma 2 of Ghosal and van der Vaart (2007) we can construct an $F'_{2,n} = \sum_{j=1}^{N_n}p_{j}\delta_{\omega_{j}}$ with $N_{n} \precsim (a_n/ \tau_{n}) \log (1/ \delta_{n})$ such that
\be \label{eq: KL3}
&&\hspace{-1cm} \norm{dF'_{1,n} \star \psi_{\sigma} \star \Normal(0,\tau_{n}^{2}) - dF'_{2,n} \star \psi_{\sigma} \star \Normal(0,\tau_{n}^{2})}_{1}  \nonumber\\
&&\precsim \norm{dF'_{1,n} \star \Normal(0,\tau_{n}^{2}) - dF'_{2,n} \star \Normal(0,\tau_{n}^{2}) }_{1} \precsim \delta_{n} \{\log (1/\delta_{n}) \}^{1/2}.
\ee

\noindent {\bf Fourth Term:} 
%
Consider disjoint balls $\{V_{j}\}_{j=1}^{N_n}$ with centers $\{\omega_{j}\}_{j=1}^{N_n}$ and diameter $\tau_{n}\delta_{n}$. 
Define $V_{0} = \rR - \cup_{j=1}^{N_n}V_{j}$. 
Extend $\{V_{j}\}_{j=1}^{N_n}$ to a partition $\{V_{j}\}_{j=1}^{M_n}$ of $[-a_{n}-1,a_{n}+1]$ such that $diam(V_{j}) = b_{2}\tau_{n}$ for all $j=N_{n}+1,\ldots,M_{n}$ for some $b_{2} \leq 1/2$. 
Define $p_{j} = 0$ for $j=N_{n}+1,\ldots, M_{n}$. 
Then applying Lemma 5 from Ghosal and van der Vaart (2007), we have
\be \label{eq: KL4}
&& \hspace{-1cm} \norm{ dF_{\alpha} \star \psi_{\sigma} \star \Normal(0,\tau_{n}^{2}) - dF'_{2,n} \star \psi_{\sigma} \star \Normal(0,\tau_{n}^{2})}_{1} \\
&& \precsim \norm{ dF_{\alpha} \star \Normal(0,\tau_{n}^{2}) - dF'_{2,n} \star \Normal(0,\tau_{n}^{2})}_{1}  \nonumber\\
&& \precsim \delta_{n} + \sum_{j=1}^{N_{n}} \abs{F_{\alpha}(V_j) - p_j}    \leq \delta_{n} + \sum_{j=1}^{N_{n}} \abs{F_{\alpha}(V_j) - p_j} + \sum_{j=N_{n}+1}^{M_{n}}  \abs{F_{\alpha}(V_j)} \nonumber\\
&&= \delta_{n} + \sum_{j=1}^{M_{n}} \abs{F_{\alpha}(V_j) - p_j}.
\ee

An upper bound to the second term in (\ref{eq: KL4}) is obtained using Lemma 10 from Ghosal and van der Vaart (2007). 
%
Let $N_n = C_{1} a_{n}\log(1/\delta_{n})/\tau_{n}$ for some $C_{1}$.
For $n$ sufficiently large and for some $b_{3} > (\gamma_{3}/c_{2}+1)/b_{1}$, we have $\delta_{n}^{b_{3}} \leq \alpha(V_{j}) \leq 1$ for all $j=1,\ldots,M_{n}$. 
Also
\bse
(M_{n}-N_{n}) &=& (2a_{n} + 2 - N_{n}\tau_{n}\delta_{n})/(b_{2}\tau_{n}) = (2a_{n}+2)/(b_{2}\tau_{n}) - N_{n}\delta_{n}/b_{2} \\
\Rightarrow  M_{n} &\precsim& N_{n} + \frac{a_{n}}{\tau_{n}} \precsim \frac{a_{n}}{\tau_{n}}\log(1/\delta_{n}) + \frac{a_{n}}{\tau_{n}} ~\precsim~ \frac{a_{n}}{\tau_{n}} \log(1/\delta_{n}).  \\
\Rightarrow  M_{n} \delta_{n}^{b_{3}} &\precsim& a_{n} \tau_{n}^{-1} \delta_{n}^{b_{3}}\log(1/\delta_{n})  = a_{0} \tau_{n}^{-\gamma_{3}/c_{2}-1} \delta_{n}^{b_{3}}\log(1/\delta_{n}) \\
&\precsim& \epsilon_{n}^{2b_{1}b_{3}-2\gamma_{3}/c_{2}-2} \{\log(1/\epsilon_{n})\}^{3+2\gamma_{3}/c_{2}} \leq 1.
\ese
According to $\Pi$, $(F_{\alpha}(V_1), \ldots, F_{\alpha}(V_{M_n})) \sim \mbox{Dir}(\alpha(V_1), \ldots, \alpha(V_{M_n}))$. Therefore, we have, using Lemma 10 from Ghosal and van der Vaart (2007),
\bse
\Pi \bigg( \sum_{j=1}^{M_{n}} \abs{ F_{\alpha}(V_j) - p_j} \leq 2\delta_{n}^{b_{3}}, \min_{1\leq j \leq M_{n}} F_{\alpha}(V_{j})\geq \frac{\delta_{n}^{2 b_{3}}}{2} \bigg) \geq C_{2}\exp\{-c_{4} M_{n} \log (1/\delta_{n})\}.
\ese
Therefore, using (\ref{eq: KL4}), we obtain, with $\Pi$-probability at least $C_{2}\exp\{-c_{4} M_{n} \log (1/\delta_{n})\}$,
\be \label{eq: KL5}
\norm{ dF_{\alpha} \star \psi_{\sigma} \star \Normal(0,\tau_{n}^{2}) - dF'_{2,n} \star \psi_{\sigma} \star \Normal(0,\tau_{n}^{2})}_{1}  \leq \delta_{n} +  2 \delta_{n}^{b_{3}}.
\ee

\noindent {\bf Fifth Term:} Using Young's inequality again, we have, for $h_{n} \in (\tau_{n},\tau_{n}+\tau_{n}^{2})$,
\be
&&\hspace{-1cm} \norm{dF_{\alpha} \star \psi_{\sigma} \star \Normal(0,\tau_{n}^{2}) - dF_{\alpha} \star  \psi_{\sigma} \star \Normal(0,h^{2})}_{1} \nonumber\\
&&\leq \norm{dF_{\alpha}\star \psi_{\sigma}}_{1} \norm{\Normal(0,\tau_{n}^{2}) - \Normal(0,h^{2})}_1 \precsim \frac{\abs{\tau_{n} - h}}{\tau_{n} \wedge h}  \precsim \tau_{n}.
\ee
Using Conditions \ref{cond: base measure} on the prior, 
we have 
\bse
&& \hspace{-1cm} \Pi(\tau_{n} < h < \tau_{n}+\tau_{n}^{2}) \succsim \int_{\tau_{n}}^{\tau_{n}+\tau_{n}^{2}} \exp(-K_{1}h^{-K_{2}})dh \succsim \tau_{n}^{2}\exp\{-K_{1}(\tau_{n}+\tau_{n}^{2})^{-K_{2}}\} \\
&&\succsim \tau_{n}^{2}\exp\{-K_{1}(2\tau_{n})^{-K_{2}}\} = \exp\{-2\log(1/\tau_{n})-K_{1}(2\tau_{n})^{-K_{2}}\} \\
&& = \exp[-4\log(1/\epsilon_{n})-4\log \{\log(1/\epsilon_{n})\}-K_{1}2^{-K_{2}} \epsilon_{n}^{-2K_{2}}\{\log(1/\epsilon_{n})\}^{2K_{2}}] \\
&& \succsim \exp[-K_{1}2^{-K_{2}} \epsilon_{n}^{-2K_{2}}\{\log(1/\epsilon_{n})\}^{2K_{2}+2}].
\ese
Also $M_{n}\log(1/\delta_{n}) \precsim (a_{n}/\tau_{n}) \{\log(1/\delta_{n})\}^{2} \precsim \tau_{n}^{-\gamma_{3}/c_{2}-1} \{\log(1/\delta_{n})\}^{2} \precsim \epsilon_{n}^{-2\gamma_{3}/c_{2}-2}\{\log(1/\epsilon_{n})\}^{4+2\gamma_{3}/c_{2}}$.  
Define $K^{\star} = (\gamma_{3}/c_{2}+1) \vee K_{2}$.
Then,  for some $C>0$, we have 
\bse
\Pi\left(\norm{ dF_{\alpha} \star \psi_{\sigma} \star \Normal(0,\tau_{n}^{2}) - dF'_{2,n} \star \psi_{\sigma} \star \Normal(0,\tau_{n}^{2})}_{1}  \leq \delta_{n} +  2 \delta_{n}^{b_{3}}\right) \Pi(\tau_{n} < h < \tau_{n}+\tau_{n}^{2}) \\
\succsim \exp[-C\epsilon_{n}^{-2K^{\star}}\{\log(1/\epsilon_{n})\}^{2+2K^{\star}}].
\ese
Also note that $\epsilon_{n}^{-2K^{\star}}\{\log(1/\epsilon_{n})\}^{(2+2K^{\star})} \precsim n\epsilon_{n}^{2}  \Leftrightarrow (\log n)^{(2+2K^{\star})(1-t)} \precsim n^{1-\gamma(2+2K^{\star})}$. 
For $\gamma \in (0,1/(2+2K^{\star})]$ this is true if $t > 1$.
Therefore, combining (\ref{eq: H}), (\ref{eq: KL1}), (\ref{eq: KL3}) and (\ref{eq: KL5}) we obtain, with $\Pi$-probability at least $\exp(-C n\epsilon_{n}^{2})$,
\be \label{eq: KL6}
d_{H}(f_W, f_{0W}) \precsim  \tau_{n}^{1/2} +\tau_{n}^{1+\gamma_{2}}  +  \delta_{n}^{1/2} \{\log (1/\delta_{n}) \}^{1/4} +  \delta_{n}^{1/2} + \delta_{n}^{b_{3}/2} \precsim  \tau_{n}^{1/2}.
\ee

For $n$ sufficiently large, we have $\sigma\delta_{n} \leq b_{2}\tau_{n}$. 
This ensures that for any $z\in [-a_{n}-1/2, a_{n}+1/2]$, for sufficiently large $n$, $V_{j}  \subset V_{z} = \{y: \abs{z-y}\leq \tau_{n}\}$ for some $j\in\{1,\ldots,M_{n}\}$ 
and $\min_{1\leq j \leq M_{n}}F_{\alpha}(V_{j}) \geq \delta_{n}^{2b_{3}}$.   
Therefore, for $\abs{\omega} \leq a_{n}$, $h \in (\tau_{n},\tau_{n}+\tau_{n}^{2})$ and $n$ sufficiently large
\bse
\frac{f_{W}(w)}{f_{0W}(w)}  &\geq&  M^{-1} \int \psi_{\sigma} (w-z)f_{X}(z) dz   \\
&=&  M^{-1}   \int  \psi_{\sigma} (w-z) \int \Normal(z-y\mid 0,h^{2})~ dF_{\alpha}(y) ~dz   \\
&=&  \frac{1}{(2\pi)^{1/2}hM}  \int  \psi_{\sigma} (w-z) \int \exp\{-(z-y)^{2}/(2h^{2})\}~ dF_{\alpha}(y) ~dz   \\
&\geq&  \frac{1}{2(2\pi)^{1/2}\tau_{n}M}  \int  \psi_{\sigma} (w-z) \int _{\abs{z-y}\leq \tau_{n}} \exp\{-(z-y)^{2}/(2\tau_{n}^{2})\}~ dF_{\alpha}(y) ~dz   \\
&\geq&  \frac{\exp(-1/2)}{2(2\pi)^{1/2}\tau_{n}M}  \int_{\abs{z-w}\leq 1/2}  \psi_{\sigma} (w-z) \int _{\abs{z-y}\leq \tau_{n}} dF_{\alpha}(y) ~dz \\
&\geq&  \frac{\exp(-1/2)\delta_{n}^{2 b_{3}}}{4(2\pi)^{1/2}\tau_{n}M}  \int_{\abs{z-w}\leq 1/2}  \psi_{\sigma} (w-z) ~dz  \\
&=&  \frac{\exp(-1/2)\delta_{n}^{2 b_{3}}}{4(2\pi)^{1/2}\tau_{n}M} \int_{\abs{z}\leq 1/2}  \psi_{\sigma} (z) ~dz = \frac{C_{3}\delta_{n}^{2 b_{3}}}{\tau_{n}} .
\ese

Take $\lambda_{n} = C_{3}\delta_{n}^{2 b_{3}}/\tau_{n} = C_{3} \epsilon_{n}^{4b_{1}b_{3}-2}\{\log(1/\epsilon_{n})\}^{2} \downarrow 0$. 
We have $\log (1/\lambda_{n}) = [\log C_{3} + (4b_{1} b_{3}-2)\log (1/\epsilon_{n}) -2\log \{\log(1/\epsilon_n)\}] \precsim \log (1/\epsilon_{n})$.
Also, for $n$ sufficiently large, $\lambda_{n}<e^{-1}$, that is, $\log(1/\lambda_{n}) > 1$. 
Therefore, with $\Pi$-probability at least $\exp(-C n\epsilon_{n}^{2})$, we have
\be \label{eq: KL7}
d_{H}^{2}(f_{W},f_{0W}) \log(1/\lambda_{n}) \leq d_{H}^{2}(f_{W},f_{0W}) \{\log(1/\lambda_{n})\}^{2} \precsim \tau_{n} \{\log(1/\epsilon_{n})\}^{2} \precsim \epsilon_{n}^{2}.
\ee

To use  Lemma \ref{lem: lemma 7 of GvdW}, we need to further provide upper bounds to 
$\int_{f_W/f_{0W} \leq \lambda}   f_{0W} \{\log (f_W/f_{0W})\}^2$ and  $\int_{f_W/f_{0W} \leq \lambda}  f_{0W} \log (f_W/f_{0W})$. 

Let $\min \{\sup f_{0X}, \sup \psi_{\sigma}\} = M < \infty$. 
We have $f_{0W}(w) = \int f_{0X}(x) \psi_{\sigma} (w-x) dx \leq M$. 
Therefore, for $\abs{w} \geq a_n$, we have
\bse
\frac{f_W(w)}{f_{0W}(w)} &\geq& M^{-1} \int _{\abs{z} \leq a_n} \psi_{\sigma} (w - z) f_{X}(z) dz
\succsim \abs{w}^{-c_{2}} \int_{\abs{z} \leq a_n} f_{X}(z) dz,
\ese
using $\abs{w - z} \leq  \abs{w} + \abs{z} \leq 2\abs{w}$, since $\abs{z} \leq a_{n} \leq \abs{w}$.
Now $\int _{\abs{z} \leq a_n} f_X(z) dz  \to 1\, \text{a.s.} ~as ~n \to \infty$.
Hence, for $\abs{w} \geq a_n$, with probability tending to $1$, $f_{W}(w)/f_{0W}(w) \succsim  \abs{w}^{-c_{2}} /  2$.
Define $B_{1,n} = \{w: \abs{w}\leq a_{n}, f_{W}(w)/f_{0W}(w) \leq \lambda_{n} \}$ and $B_{2,n} = \{w: \abs{w}\geq a_{n}, f_{W}(w)/f_{0W}(w) \leq \lambda_{n} \}$. 
By definition of $\lambda_{n}$, $B_{1,n}$ is actually empty. 
Using (\ref{eq: tail both}) we also have $\Pr_{f_{0W}}(B_{2,n})  \leq \Pr_{f_{0W}}(\abs{W}>a_{n}) \leq \Pr_{f_{0X}}(\abs{X}>a_{n}/2)+\Pr_{\psi_{\sigma}}(\abs{U} >a_{n}/2) \precsim \tau_{n}^{2+2\gamma_{2}}$. 
Therefore, with probability tending to $1$, we have
\be\label{eq: KL8}
&& \hspace{-1cm} \int_{\{f_{W}/f_{0W} \leq \lambda_{n}\}}  f_{0W}  \bigg(\log \frac{f_{0W}}{f_W}\bigg)^2  ~=~ 
\int_{B_{1,n}\cup B_{2,n}}  f_{0W}  \bigg(\log \frac{f_{0W}}{f_W}\bigg)^2  =  \int_{B_{2,n}}  f_{0W}  \bigg(\log \frac{f_{0W}}{f_W}\bigg)^{2}  \nonumber\\
  &\precsim&  2(\log 2)^{2}\int_{B_{2,n}} f_{0W} (w) dw  + 2c_{2}^{2} \int_{B_{2,n}}  (\log \abs{w}) f_{0W} (w) dw  \nonumber\\
  &\precsim&  \Pr_{f_{0W}}(B_{2,n}) + \Pr_{f_{0W}}^{1/2} (B_{2,n}) E_{f_{0W}} (\log\abs{W})^{2} \precsim  \tau_{n}^{1+\gamma_{2}} 
  \precsim \tau_{n}
\ee
Also, for $n$ sufficiently large, $\lambda_{n}<e^{-1}$, and we have $\log(f_{0W}/f_{W}) 1(f_{W}/f_{0W} < \lambda_{n}) \leq \{\log(f_{W}/f_{0W})\}^{2}1(f_{W}/f_{0W} < \lambda_{n})$.
Therefore, with $\Pi$-probability at least $\exp(-Cn\epsilon_{n}^{2})$,
\be \label{eq: KL9}
\int_{\{f_{W}/f_{0W} < \lambda_{n}\}} f_{0W}\left(\log\frac{f_{0W}}{f_{W}}\right) \leq \int_{\{f_{W}/f_{0W} < \lambda_{n}\}} f_{0W} \bigg(\log \frac{f_{0W}}{f_W}\bigg)^2 \precsim  \tau_{n}^{1+\gamma_{2}} \precsim \epsilon_{n}^{2}.
\ee 
Finally, combining (\ref{eq: KL7}) and (\ref{eq: KL9}), we have
$\int f_{0W}  \log (f_{0W}/f_{W}) \precsim \int  f_{0W}  \{\log (f_{0W}/f_{W})\}^2  \precsim \epsilon_{n}^{2}$ with $\Pi$-probability at least $\exp(-Cn\epsilon_{n}^{2})$. 
That is the KL conditions hold for $\epsilon_{n} = n^{-\gamma} (\log ~ n)^{t}$ provided $\gamma \in (0,1/(2+2K^{\star})]$ and $t > 1$.

\section{Proof of the Main Theorems}  \label{sec: pf main thm}

Theorem \ref{thm: main thm supersmooth} follows by verifying the sufficient conditions stated in Theorem \ref{thm: suff} 
for the respective DPMM priors as specified in Section \ref{sec: main thms}.
Condition (\ref{eq: suff2}) of Theorem \ref{thm: suff} has already been verified in the subsection \ref{sec: KL}, respectively, for the choice of $\epsilon_{n}$ in subsection \ref{sec: KL}.
Condition (\ref{eq: bias of DKE}) follows trivially choosing $h_{n}\simeq (\log n)^{-1/\beta}$.
Then it remains to be shown that (\ref{eq: estconc}) holds for $\xi_{n} \simeq (\log n)^{-\eta/\beta}$ and the choice of $\epsilon_{n}$ specified in subsection \ref{sec: KL}.\\

\noindent{\bf Case $2 \leq p <\infty$ - Application of the inequality due to Talagrand (1996):} \\
We apply (\ref{eq: tallagrand}) to verify (\ref{eq: estconc}).
Using $(a+b)^{1/2} \leq a^{1/2} + b^{1/2}$ and $(ab)^{1/2} \leq (a+b)/2$, we have
\bse
&& \hspace{-1cm} (2K_{2}x)^{1/2} = [2 \{n K_{3} + 2 K_{1}n E\norm{f_{n, X} - E_{f_{W}} f_{n,X}}_{p}\}x]^{1/2} \\
&& \leq  (2n K_{3} x)^{1/2} + \{4 K_{1} x n E\norm{f_{n, X}(x) - E_{f_{W}} f_{n,X}(x)}_{p}\}^{1/2} \\
&& \leq (2n K_{3} x)^{1/2} + 2 K_{1} x + n E\norm{f_{n, X} - E_{f_{W}} f_{n,X}}_{p}/2.
\ese
Therefore,
\bse
&&\hspace{-.1cm} Pr_{f_{W}}\{ n\norm{f_{n,X} - E_{f_{W}} f_{n,X}}_{p} \geq (3/2) n E_{f_{W}} \norm{f_{n,X} - E_{f_{W}} f_{n,X}}_{p} + (2nK_{3}x)^{1/2} + (7/3)K_{1}x\}   \\
&& \leq Pr_{f_{W}}\{ n\norm{f_{n,X} - E_{f_{W}} f_{n,X}}_{p} \geq n E_{f_{W}} \norm{f_{n,X} - E_{f_{W}} f_{n,X}}_{p} + (2K_{2}x)^{1/2} + K_{1}x/3\}   \leq 2e^{-x}.
\ese
Take $x = Ln\epsilon_{n}^{2}$ with $L>(4+C)$, $\xi_{n} = (\log n)^{-\eta/\beta}$ and $\epsilon_{n}$ as in subsection \ref{sec: KL}. 
Take $h_{n} = \{2/(\gamma\varrho)\}^{1/\beta} (\log n)^{-1/\beta}$.
Then, it is easy to verify that 
\bse
nD_{1}\xi_{n} \succsim (3/2) n^{1/2}h_{n}^{\varpi_{p}(\beta_{0})} \exp(h_{n}^{-\beta}/\varrho) + \{2 n L n \epsilon_{n}^{2} h_{n}^{2\varpi_{\infty}(\beta_{0})}\exp(2h_{n}^{-\beta}/\varrho)\}^{1/2} \\ ~~~~~~ + (7/3)h_{n}^{\varpi_{\infty}(\beta_{0})} \exp(h_{n}^{-\beta}/\varrho)Ln\epsilon_{n}^{2}.
\ese

\noindent{\bf Case $p=\infty$ - Application of the inequality due to Dvoretzsky, et al. (1956):} \\ 
In this case, we can apply (\ref{eq: Dvoretzsky}) to verify (\ref{eq: estconc}).
Take $h_{n} = \{4/(\gamma\varrho)\}^{1/\beta}(\log n)^{-1/\beta}$.
Then 
\bse
n\xi_{n}^{2}h_{n}^{2}/V_{n}^{2} ~ \succsim ~ n\xi_{n}^{2}h_{n}^{-\varpi^{\infty}(\beta_{0}) }\exp(-4h_{n}^{-\beta}/\varrho) ~ \succsim ~ n (\log n)^{\varpi^{\infty}(\beta_{0})/\beta-2\eta/\beta} n^{-\gamma} ~ \succsim ~ n^{1-2\gamma} (\log n)^{2t} = n\epsilon_{n}^{2}.
\ese
Combining the different cases, Theorems \ref{thm: main thm supersmooth} follows.

\section{Accelerated Rates} \label{sec: acc rates}

The optimal rate of convergence of deconvolution estimators for supersmooth errors is extremely slow.
It is, therefore, desirable to determine how small $\sigma$ should be for deconvolution to be practically feasible 
and for the deconvolution estimator to converge as fast as an ordinary density estimator (Fan, 1992).
The following theorem shows that if $\sigma$ decreases at a certain rate as a function of the sample size, 
then deconvolution is as difficult as Bayesian density estimation in terms of the rate of convergence.

\begin{Thm} \label{thm: density estimation rates}
Let $2 \leq p < \infty$, $f_{0X} \in \C_{p}^{\eta}(\rR) \cap \C_{\infty}^{\eta}(\rR)$ and $\Pi$ be a DPMM prior with a $\Normal(\mu_{0},\sigma_{0}^{2})$ base measure and an $\Exp(\lambda)$ prior on $h$. 
Assume there exist positive constants $c_{1}, c_{2}, c_{3}$ and $T$ such that 
\be \label{eq: exponential tails} 
f_{0X}(x) \leq c_{1} \exp(-c_{2}\abs{x}^{c_{3}}), \quad \abs{x} \geq T,~~~f \in \{f_{0X},\psi_{\sigma_{n}}\},
\ee
that is, both $f_{0X}$ and $\psi_{\sigma_{n}}$ have exponentially decaying tails.
Let $\sigma = \sigma_{n} \simeq n^{-1/(2\eta+1)} (\log n)^{t/\eta}$ where $t > (2+1/\eta+1/c_{3})/(2+1/\eta)$.
Also let at least one $f$ in $\{f_{0X},\psi_{\sigma_{n}}\}$ satisfy 
\be \label{cond on true1}
\int \left\{\abs{f^{(k)}}/f\right\}^{(2\eta+\rho)/k} f < \infty, ~~~ k\in \{1,2,\ldots,\lfloor\eta\rfloor\}, ~~~~~~\int f^{1-(2\eta+\rho)/\eta} < \infty,
\ee
for some $\rho > 0$. 
Then there exists an $M \in (0,\infty)$ such that
\be \label{eq: thm density estimation rates}
E_{f_{0W}}\Pi\{ \norm{f_X - f_{0X}}_p \geq M n^{-\eta/(2\eta+1)} (\log n)^{t} \mid \bW_{1:n}\} \to 0.
\ee
\end{Thm}
\MyProof
A proof is obtained by verifying the sufficient conditions stated in Theorem \ref{thm: suff} for $\xi_{n} = n^{-\eta/(2\eta+1)} (\log n)^{t}$.
Let $h_{n} \simeq \sigma_{n}$. 
Then (\ref{eq: bias of DKE}) follows by the calculations of Section \ref{sec: DKE}.
The variance bounds obtained in Section \ref{sec: concentration bounds for DKE} can now be refined as follows.  
Using the fact that $\phi_{\sigma_{n}}(t/h_{n}) = \phi(t\sigma_{n}/h_{n})$, 
we now have $\sup_{k\in \K_{p}}\abs{k} \precsim K_{1} = O(1)$; $\sup_{k\in \K_{p}}E_{f_{W}}k^{2}(W) \precsim K_{2} = O(1)$;
and $E_{f_{W}}\norm{n\wh{F}_{n,W}}_{\K_{p}} \precsim n^{1/2}$.
Using Lemma B3 from Shen, et al. (2013) and by minor adjustment of the proof of Theorem 4 from the same paper we can verify 
the KL condition (\ref{eq: suff2}) with $\epsilon_{n} = n^{-\eta/(2\eta+1)} (\log n)^{t}$.
Condition (\ref{eq: estconc}) is then readily verified by application of Talagrand's inequality since $n\xi_{n} \succsim n^{1/2} + n\epsilon_{n} + n\epsilon_{n}^{2}$.

\section{Discussion}\label{sec: discussion}

In this article we provided a set of sufficient conditions for the posterior of a Bayesian deconvolution estimator to converge to the true density of interest in $L_{p}$ norm, 
assuming the measurement error density to be known and supersmooth.
We showed that under a minimal polynomial tails assumptions on the density of interest and the error density and for $2 \leq p \leq \infty$, 
the sufficient conditions hold for a location mixture of Normals prior on the density of interest induced by a Dirichlet process with a conjugate Normal base measure and a conjugate inverse-gamma prior on the bandwidth parameter.
The posterior of such a DPMM prior actually attains the minimax optimal rate of convergence $(\log n)^{-\eta/\beta}$, where $\eta$ is the smoothness of the true density and $\beta$ is the degree of smoothness of the error distribution. 
The minimax rate is, in fact, achieved adaptively - no prior knowledge of $\eta$ is necessary. 

The case of ordinary smooth errors, where the tails of the characteristic function of the error density decay polynomially, is not considered in this article. 
The proof of convergence results for ordinary smooth errors would require much stronger prior concentration bounds for Kullback-Leibler balls around the true convolved density. 
We are pursuing this problem as the subject of separate research.

\section*{Acknowledgments}
The work is part of the first author's PhD dissertation.
Carroll's research was supported in part by a grant R37-CA057030 from the National Cancer Institute.
Mallick's research was supported in part by NSF grant DMS0914951.
This publication is based in part on work supported by Award Number KUS-CI-016-04, made by King Abdullah University of Science and Technology (KAUST).


\section*{References}
\refmark
Bousquet, O. (2003). Concentration inequalities for sub-additive functions using the entropy method. 
\emph{Progress in Probability}, 56, 213-247.
\refmark
Carroll, R. J. and Hall, P. (1988). Optimal rates of convergence for deconvolving a density.
\JASA, 83, 1184-1186.
\refmark
Dvoretzsky, A.,  Kiefer, J. and Wolfowitz, J. (1956). Asymptotic minimax character of the sample distribution function and of the classical multinomial estimator.
\emph{Annals of Mathematical Statistics}, 27, 642-669.
\refmark
Escobar, M. D. and West, M. (1995). Bayesian density estimation and inference using mixtures.
\JASA, 90, 577-588.
\refmark
Fan, J. (1988). Optimal global rates of convergence for nonparametric deconvolution problems. 
emph{Technical Report 162}, Department of Statistics, University of California, Berkeley.
\refmark
Fan, J. (1991a). Global behavior of deconvolution kernel estimates. 
\SSNC, 1, 541-551.
\refmark
Fan, J. (1991b). On the optimal rates of convergence for nonparametric deconvolution problems. 
\ANNALS, 19, 1257-1272.
\refmark
Fan, J. (1992). Deconvolution with supersmooth errors. 
\CANADAJS, 20, 155-169.
\refmark
Fan, J. and Hu (1992). Bias correction and higher order kernel functions. 
\STATL, 13, 235-243.
\refmark
Ferguson, T. F. (1973). A Bayesian analysis of some nonparametric problems.
\ANNALS, 1, 209-230.
\refmark
Ghosal, S., Ghosh, J. K. and Ramamoorthi, R. V. (2000). Posterior consistency of Dirichlet mixtures in density estimation. 
\ANNALS, 27, 143-158.
\refmark
Ghosal, S., Ghosh, J. K. and Van Der Vaart, A. W. (2000). Convergence rates of posterior distributions. 
\ANNALS, 28, 500-531.
\refmark
Ghosal, S. and van der Vaart, A. W. (2007). Posterior convergence rates of Dirichlet mixtures at smooth densities. 
\ANNALS, 35, 697-723.
\refmark
Gin\'{e}, E. and Mason, D. M. (2007). On local U-statistic processes and the estimation of densities of functions of several sample variables. 
\ANNALS, 35, 1105-1145.
\refmark
Gin\'{e}, E. and Nickl, R. (2012). Rates of contraction for posterior distributions in $L_{r}$-metrics, $1\leq r \leq \infty$. 
\ANNALS, 39, 2883-2911.
\refmark
Knapik, W., van der Vaart, A. W. and van Zanten, J. H. (2011). Bayesian inverse problems with Gaussian priors. 
\ANNALS, 39, 2626-2657.
\refmark
Kruijer, W., Rousseau, J. and van der Vaart, A. W. (2010). Adaptive Bayesian density estimation with location-scale mixtures. 
\EJS, 4, 1225-1257.
\refmark
Lo, A. Y. (1984). On a class of Bayesian nonparametric estimates. I: Density estimates.
\ANNALS, 12, 351-357.
\refmark
Prakasa Rao, B. L. S. (1983). \emph{Nonparametric Functional Estimation}.
Academic, New York.
\refmark
Ray, K. (2012). Bayesian inverse problems with non-conjugate priors. arXiv preprint arXiv:1209.6156.
\refmark
Sarkar, A., Mallick, B. K., Staudenmayer, J., Pati, D. and Carroll, R. J. (2013). Bayesian semiparametric density deconvolution in the presence of conditionally heteroscedastic measurement errors.
Unpublished manuscript.
\refmark
Shen, W., Tokdar, S. T. and Ghosal, S. (2013). Adaptive Bayesian multivariate density estimation with Dirichlet mixtures. 
\BIOK, 1-18.
\refmark
Shen, W. and Wasserman, L. (2001). Rates of convergence of posterior distribution.
\ANNALS, 29, 687-714.
\refmark
Stefanski, L. A. and Carroll, R. J. (1990). Deconvoluting kernel density estimators.
{\it Statistics}, 21, 169-184.
\refmark
Talagrand, M. (1996). New concentration inequalities in product spaces.
\emph{Inventiones Mathematicae}, 126, 505-563.

\newcommand{\Appendix}{\appendix\def\thesection{Appendix~\Alph{section}}\def\thesubsection{\Alph{section}.\arabic{subsection}}}
\section*{Appendix}
\begin{appendix}
\Appendix
\renewcommand{\theequation}{A.\arabic{equation}}
\setcounter{equation}{0}
\baselineskip=18pt

\section{Proof of Theorem \ref{thm: suff}} \label{pf: main thm suff}
We will use the proofs of Theorem 2.1 in Ghosal, et al. (2000) in conjunction with the proof of Theorem 3.1 in Ray (2012). 
Define $\S_{n} = \{f_W: f_W = f_X \star \psi_{\sigma}, \norm{f_X-f_{0X}}_{p}\geq M\xi_{n} \}$, where $M$ is a positive constant to be chosen later. 
A key step is to construct tests (indicator functions) $\Phi_n$ as a function of the observed data $\bW_{1:n}$ such that 
\bse \label{eq:tests}
E_{f_{0W}} (\Phi_n) \to 0,  
\quad \sup_{ f_W\in \calP_{n}\cap\S_{n}}  E_{f_W} (1 - \Phi_n) \leq \exp\{-(4+C)n\epsilon_{n}^2\} .
 \ese
Consider the test $\Phi_n = 1 _{\norm{f_{n, X} - f_{0X}}_p > M_{1} \xi_n}$, where $M_{1}$ is a constant to be chosen later.
\bse
E_{f_{0W}} \Phi_{n} &=& \Pr_{f_{0W}}\left(\norm{f_{n, X} - f_{0X}}_p \geq M_{1} \xi_n\right) \\
&\leq&  \Pr_{f_{0W}}\left(\norm{f_{n, X} - E_{f_{0W}} f_{n,X}}_p + \norm{E_{f_{0W}} f_{n,X}-f_{0X}}_p \geq M_{1} \xi_n\right) \\
&=&  \Pr_{f_{0W}}\left(\norm{f_{n, X} - E_{f_{0W}} f_{n,X}}_p \geq M_{1} \xi_n - \norm{E_{f_{0W}} f_{n,X}-f_{0X}}_p \right) \\
&\leq&  \Pr_{f_{0W}}\left\{\norm{f_{n, X} - E_{f_{0W}} f_{n,X}}_p \geq (M_{1}-D_{0}) \xi_n \right\}.
\ese
Taking $M_{1} \geq (D_0+D_1)$ we have
\be
E_{f_{0W}} \Phi_{n} \leq \exp\{-(4+C)n\epsilon_{n}^2\}. \label{eq:suffpf1}
\ee
We have $\norm{f_{X}-f_{0X}}_p \leq \norm{f_{0X}-f_{n,X}}_{p} + \norm{f_{X}-E_{f_{W}}f_{n,X}}_{p} + \norm{E_{f_{W}}f_{n,X} - f_{n,X}}_{p}$.
For any $f_{W}\in \calP_{n}\cap \S_{n}$,
\bse
&&\hspace{-1cm} E_{f_{W}} (1-\Phi_{n})  = \Pr_{f_{W}}\left(\norm{f_{n,X}-f_{0X}}_{p} \leq M_{1} \xi_{n} \right) \\
&\leq& \Pr_{f_{W}} \left(\norm{f_{X}-f_{0X}}_p - \norm{f_{X}-E_{f_{W}}f_{n,X}}_{p} - \norm{E_{f_{W}}f_{n,X} - f_{n,X}}_{p} \leq M_{1} \xi_{n} \right) \\
&\leq& \Pr_{f_{W}} \left(M\xi_{n} - D\xi_{n} - \norm{E_{f_{W}}f_{n,X} - f_{n,X}}_{p} \leq M_{1} \xi_{n} \right)  \\
&=& \Pr_{f_{W}} \left\{\norm{E_{f_{W}}f_{n,X} - f_{n,X}}_{p} \geq (M - D - M_{1}) \xi_{n} \right\}.
\ese
Taking $M \geq (D+D_1+M_1)$ we have
\be
\sup_{f_{W}\in\calP_{n}\cap \S_{n}}E_{f_{W}} (1-\Phi_{n}) \leq \exp\{-(C+4)n\epsilon_{n}^2\}. \label{eq:suffpf2}
\ee
Combining the hypothesis (\ref{eq: suff1}) with (\ref{eq:suffpf2}) we have
\bse
&&\hspace{-1cm} E_{f_{0W}} \int_{\S_{n}} \prod_{i=1}^{n}\{f_{W}(W_{i})/f_{0W}(W_{i})\}d\Pi_{n}(f_W)(1-\Phi_{n})  \\
&\leq& \Pi_{n}(\calP_{n}^{c}) + \int_{\calP_{n}\cap \S_{n}} (1-\Phi_{n}) \prod_{i=1}^{n} f_{W}(W_{i}) d\mu ~ d\Pi_{n}(f_W)  \\
&\leq& \Pi_{n}(\calP_{n}^{c}) + \int_{\calP_{n}\cap \S_{n}} E_{f_{W}}(1-\Phi_{n}) ~ d\Pi_{n}(f_W)  \\
&\leq& 2\exp\{-n\epsilon_{n}^{2}(4+C)\}.
\ese
Denote $B_n = \{f_W: \int f_{0W} \log (f_{0W}/f_{W}) < \epsilon_{n}^2,  \int f_{0W} \{\log (f_{0W}/f_{W})\}^2 < \epsilon_{n}^2\}$.
Using hypothesis (\ref{eq: suff2}) and Lemma 8.1 of Ghosal, et al. (2000), we have, with probability tending to 1,
\bse
&&\hspace{-1cm} \int \prod_{i=1}^{n}\{f_{W}(W_i)/f_{0W}(W_i)\} d\Pi_{n}(f_{W}) \geq \int_{B_n}  \prod_{i=1}^{n}\{f_{W}(W_i)/f_{0W}(W_i)\} d\Pi_{n}(f_{W}) \\
&=&  \Pi_{n}(B_n) \int_{B_n} \prod_{i=1}^{n}\{f_{W}(W_i)/f_{0W}(W_i)\} d\Pi_{n}(f_{W})/ \Pi_{n}(B_n)   \\
&\geq&  \exp(-nC\epsilon_{n}^{2}) \exp(-2n\epsilon_{n}^{2}) = \exp\{-n\epsilon_{n}^2(2+C)\}.
\ese
\bse
&&\hspace{-1cm} E_{f_{0W}} \Pi_{n}(\norm{f_{X}-f_{0X}}_{p} \geq M\epsilon_{n}\mid \bW_{1:n}) (1-\Phi_{n})1_{A_n}  \\
&=& \int_{A_n} \Pi_{n}(\norm{f_{X}-f_{0X}}_{p} \geq M\epsilon_{n}\mid \bW_{1:n})(1-\Phi_n)  \prod_{i=1}^{n}f_{0W}(W_i)d\mu  \\
&=& \int_{A_n}     \frac{\int_{\S_{n}}\prod_{i=1}^{n}\{f_{W}(W_i)/f_{0W}(W_i)\} d\Pi_{n}(f_{W})}{\int \prod_{i=1}^{n}\{f_{W}(W_i)/f_{0W}(W_i)\} d\Pi_{n}(f_{W})}         (1-\Phi_n)         \prod_{i=1}^{n}f_{0W}(W_i)d\mu  \\
&\leq&  \exp\{n\epsilon_{n}^{2}(2+C)\}    \int_{A_n}   \int_{\S_{n}}\prod_{i=1}^{n}\{f_{W}(W_i)/f_{0W}(W_i)\} d\Pi_{n}(f_{W})   (1-\Phi_n) \prod_{i=1}^{n}f_{0W}(W_i)d\mu  \\
&\leq&  \exp\{n\epsilon_{n}^{2}(2+C)\}    ~   E_{f_{0W}}  \int_{\S_{n}}\prod_{i=1}^{n}\{f_{W}(W_i)/f_{0W}(W_i)\} d\Pi_{n}(f_{W})   (1-\Phi_n)   \\
&\leq&  \exp\{n\epsilon_{n}^{2}(2+C)\}    ~   2 \exp\{-n\epsilon_{n}^{2}(4+C)\}  = 2\exp(-2n\epsilon_{n}^{2}) \to 0.
\ese

\end{appendix}

\end{document}